\documentclass[11pt,final,3p]{elsarticle}
\usepackage{amssymb}
\usepackage{amsthm}
\usepackage{amsmath}
\usepackage{latexsym}
\usepackage{amsmath}
\usepackage{xcolor}
\usepackage{graphicx}
\usepackage{indentfirst}
\usepackage[OT2,OT1]{fontenc}
\usepackage{mathrsfs}
\biboptions{numbers,sort&compress}
\usepackage{pdfsync}
\usepackage{hyperref}
\hypersetup{hypertex=true,
	colorlinks=true,
	linkcolor=blue,
	anchorcolor=blue,
	citecolor=blue}
\allowdisplaybreaks
\usepackage[utf8]{inputenc}

\renewcommand{\thefootnote}{\arabic{footnote}}

\numberwithin{equation}{section}

\newcommand\blfootnote[1]{%
	\begingroup
	\renewcommand\thefootnote{}\footnote{#1}%
	\addtocounter{footnote}{-1}%
	\endgroup
}
\journal{?}
\begin{document}
\begin{frontmatter}
\title{ Stochastic functional partial differential equations with  monotone coefficients: Poisson stability measures, exponential mixing and limit theorems }

\author{{ \blfootnote{$^{*}$Corresponding author } Shuaishuai Lu$^{a}$ \footnote{ E-mail address : luss23@mails.jlu.edu.cn}
		,~ Xue Yang$^{a,*}$}  \footnote{E-mail address : xueyang@jlu.edu.cn},~ Yong Li$^{a,b}$  \footnote{E-mail address : liyong@jlu.edu.cn}\\
	{$^{a}$College of Mathematics, Jilin University,} {Changchun 130012, P. R. China.}\\
	{$^{b}$School of Mathematics and Statistics, and Center for Mathematics and Interdisciplinary Sciences, Northeast Normal University,}
	{Changchun 130024, P. R. China.}
}

\begin{abstract}
This paper examines Poisson stable (including stationary, periodic,  almost periodic, Levitan almost periodic, Bohr almost automorphic, pseudo-periodic,  Birkhoff recurrent, pseudo-recurrent, etc.) measures and limit theorems  for stochastic functional partial differential equations(SFPDEs) with  monotone coefficients. We first show the existence and uniqueness of   entrance measure $\mu _{t}$ for SFPDEs  by dissipative method  (or remoting start). Then, with the help of Shcherbakov’s comparability method in character of recurrence, we prove  that the entrance measure inherits the same  recurrence of coefficients. Thirdly, we show the tightness of the set of measures $\mu _{t}$. As a result, any sequence of the average of  $\{\mu _{t}\}_{t\in\mathbb{R} }$ have the limit point $\mu ^{*}$. Further, we study the uniform exponential mixing of the measure $\mu ^{*}$ in the sense of Wasserstein metric. Fourthly, under  uniform exponential mixing and Markov property, we establish the strong law of large numbers, the central limit theorem and estimate the corresponding rates of convergence for solution maps  of SFPDEs.  Finally, we give  applications of stochastic generalized porous media equations with delay to illustrate of our results.
~\\
~\\
 \textbf{keywords}:  SFPDEs with  monotone coefficients; Poisson stable measures; Exponential mixing; Strong law of large numbers; Central limit theorem.
\end{abstract}
\end{frontmatter}
\section{\textup{Introduction}}
	
 In the present paper, we consider the following stochastic functional partial differential equation with  monotone coefficients:
\begin{equation}\label{g39}
\begin{cases}
 \text{d}u(t)=[A(t,u)+f(t,u_{t})]\text{d}t+g(t,u_{t})\text{d}W(t), \\
u_{s} =\varphi\in \mathcal{H},
\end{cases}
	\end{equation}
where $A$ satisfies some monotone conditions and $W(t)(t \in \mathbb{R} )$ is a two-sided cylindrical $Q$-Wiener process with $Q = I$ on a separable Hilbert space $(K,\langle \cdot , \cdot \rangle_{ K })$. Now if  $u(t;s,\varphi)$ is the solution with the initial segment $\varphi$ at $s$ of \eqref{g39}, we have a solution map $u_{t}(s, \varphi)$ by setting  $u_{t} =\{u_{t} (\theta )\}=u(t+\theta )(-\tau \le \theta \le 0)$.

Recurrent motion is an important research topic in dynamical systems, which plays an important role in the discussion of stability.  For nondeterministic systems, Kolmogorov first put forward the corresponding concept of recurrence in 1930s.  However, so far there has been a few of literature about recurrence with well properties such as  periodicity and almost periodicity, besides stationarity. In addition, it is  difficult to study the  more general recurrence, such as Levitan almost periodic, Bohr almost automorphic, pseudo-periodic, Birkhoff recurrent, pseudo-recurrent, Poisson stable, etc. It, therefore, needs to  establish a unified framework to study recurrence, particularly to determine recurrence  according to the recurrence of the coefficients for stochastic systems. In this paper, we will observe  this question and study the recurrence property of system \eqref{g39} in the perspective of measures.

It is worth noting that the solution of \eqref{g39} is non-Markov, so we have to consider the solution map with Markov properties alternatively. Then, the measure with the same recurrence property as the coefficients is  analyzed. In addition,  establishing  the limit theorems of Markov process is one of the central themes of probability theory, especially the strong law of large numbers(SLLN) and the central limit theorem(CLT), which describe the long-term behavior of stochastic processes. Therefore, in this paper, we discuss a limit measure with exponential mixing property by studying the compactness of the set of  measures, and aim at  establishing SLLN and CLT of solution maps.

 Let us make a further recall. Poisson stable solutions of system was used to express the most general recurrence, which was first introduced by Poincar$\acute{\text{e}} $ at the end of 19th century. After that, the properties of Poisson stable functions have been continuously improved, for example, see \cite{R221,R222,R223,R225,R226,R227}. Ones have considered the following classes of  Poisson stable solutions for stochastic differential equations: periodic \cite{R12,R13,R25,R26,R27}, quasi-periodic \cite{R1}, almost periodic \cite{R16,R17}, almost automorphic \cite{R18,R23} and Poisson stable \cite{R69,R3,R22,R28}, among others. It is worth noting that the above articles mainly focus on studying the recurrence property of solutions. According to the relationship between the law of solutions and measures of stochastic differential system, there are some researches around various measures of the system. Let us review some related literature in recurrent measures.  The papers \cite{R12,R31,R1,R47} studied the existence of periodic measures and quasi-periodic measures of finite-dimensional stochastic equations. For periodic measures of stochastic partial differential equations, it follows from \cite{R33,R35,R36}. However, for stochastic partial differential, there are very few publications about measures which have the more general recurrence properties.

In many practical applications, such as physical models, population prediction models and so on, the evolution of the systems is not only influenced by the present state, but also related to some of its past states. In other words, delays are ubiquitous in these systems. In addition, many phenomena in real life are inherently random. There are some studies in the stochastic dynamical properties of  functional differential equations with time delay  in recent years, see, for example \cite{R37,R39,R51,R42,R43,R45,R46,R48,R49,R111} et al. Simultaneously, some researches on  recurrence property for stochastic functional differential equations(SFDEs) have been underway, see, e.g., \cite{R50,R33,R35,R36,R52,R53}.  But there is no unified framework to study Poisson  stable measures for SFDEs. Motivated by the work of \cite{R69}, we try to develop Shcherbakov’s ideas and methods to study Poisson  stable measures of  \eqref{g39}. More precisely, we will establish suitable conditions to ensure the existence of a measure with the same recurrence property as the coefficients.

In this paper, we first prove that there exists an entrance measure $\mu _{t} (t \in \mathbb{R} )$ of \eqref{g39} by dissipative method, according to the concept of entrance measures for stochastic differential equations introduced by \cite{R1}. Then we show that the entrance measure $\mu _{t} $ inherits the same recurrent  properties as coefficients  by Shcherbakov's comparability method in the character of recurrence. In order to  obtain  the recurrent measure $\mu _{t} $ of \eqref{g39}, there are two main problems to be solved. The first is the construction of entrance measure.  Indeed, the entrance measure $\mu _{t} $ we obtained is the law of a special class of solution maps, i.e, $\mu _{t}(\Gamma )=\mathbb{P} (\omega :\mathcal{U}_{t}\in\Gamma )$ where $u_{t}(s,\varphi )\overset{\mathcal{L}  ^{2}}{\rightarrow} \mathcal{U}_{t} $ as $s\to -\infty $. The other of the  essential difficulty  is to prove   the weak compactness of the set of all entrance measures. In \cite{R3}, the tightness of the law of solutions was proved by the compact imbedding. On this basis, inspired by \cite{R33}, we obtain the  tightness of the law of  solution maps (i.e., the weak compactness of the set of all
entrance measures) by Krylov-Bogolyubov’s method and Arzelà-Ascoli’s theorem. We thus obtain  the first major result of this paper (see Theorem 4.5 and  Corollary 4.6 for details), which asserts the recurrent entrance measure.

When the system \eqref{g39} is homogeneous, i.e., $A$, $f$ and $g$ do not depend on $t$, in this paper, another deep aims are to establish SLLN and CLT for SFPDEs with monotone coefficients. Since Doeblin established the limit theorems of continuous and discrete Markov processes in \cite{R55}, various problems about SLLN and CLT have been greatly developed, see \cite{R56,R57,R58,R59,R60,R61,R7,R331,R332,R333} et al. It is worth noting that the above results are based on the fact that the solution of the system is a Markov process with exponential mixing. Obviously, those conclusions can not be applied to \eqref{g39}, because the solution of stochastic functional system with  delay depends on the history,  thus it is non-Markov, consequently we have to study  the solution map alternatively. At present, there are a few results for SFDEs with delay, see \cite{R8,R62,R65}. We note that the conclusions in these publications,  only the limit theorem of finite-dimensional SFDEs defined on $\mathbb{R}^{n}$  was considered. As far as we know, there is no conclusion about limit theorems of stochastic partial differential equations with monotone coefficients. Thus, in this paper, we will establish SLLN and CLT for solution map of  \eqref{g39}.

 As mentioned above, the exponential mixing plays an important role in analyzing SLLN and CLT of system \eqref{g39}. The conclusions about exponential mixing of Markov processes can be found in \cite{R33,R61,R67,R68,R66,R112,R113,R555}, etc. In order to obtain exponential mixing, under  some dissipative condition, we will investigate the large-time behavior of the measure $	\mu ^{L}=\frac{1}{L} \int_{0}^{L} \mu _{t}\text{d}t$.
It has been noted that the tightness of $\{\mu _{t}\}_{t\in \mathbb{R}}$ implies that $$\{\mu ^{L}=\frac{1}{L} \int_{0}^{L} \mu _{t}\text{d}t:L\in \mathbb{N}^{+} \}$$
 is tight. Then we obtain that when $L\to \infty $, there exists a weak limit $\mu^{*}$ which is an invariant measure and satisfies uniform exponentially mixing in the sense of Wasserstein metric(see Theorem 5.2 for details).

 Finally, under uniform exponential mixing and Markov properties of the solution map, we establish the strong law of large numbers, the central limit theorem and the corresponding rates of convergence(see Theorem 6.1 and Theorem 6.6 for details):
 \begin{enumerate}[(1)]
		\item Strong law of large numbers:
$$\frac{1}{t}\int_{0}^{t}F(u_{s}(\varphi ))\text{d}s \to \int _{\mathcal{H} } F(\phi ) \mu ^{*}(\text{d}\phi ) \quad \text{as}\quad t\to \infty ,\quad \mathbb{P} -a.s.,$$
where $F$ is the observation function;
        \item  Central limit theorem:
        $$\frac{1}{\sqrt{t} }\int_{0}^{t}[F(u_{r}(\varphi ))-\int _{\mathcal{H} } F(\phi ) \mu ^{*}(\text{d}\phi )]\text{d}r \overset{W}{\to}\Pi ,$$
        where $\overset{W}{\to}$ means weak convergence and $\Pi $ is a normal random variable.
	\end{enumerate}

  The structure of the paper is as follows.  In the next section, we introduce some definitions, notation, lemmas and  the basic concepts   of  Poisson stable functions, as well as  comparable (uniformly comparable) methods of Shcherbakov. In Section 3, we examine the existence and uniqueness
of the entrance measure for  SFPEDs with  monotone coefficients. In Section 4, we show that the entrance measure is uniformly compatible with coefficients, so it possesses the same character of recurrence as the coefficients.  In Section 5, we discuss the exponential mixing of the limit measure. In Section 6, we establish the SLLN and CLT for  autonomous systems. In Section 7, we give an application in stochastic generalized porous media equations.
	\section{\textup{Preliminaries}}
	
	Let $(U,\left \| \cdot  \right \| _{U})$, $(K,\left \| \cdot  \right \| _{K})$  be separable Hilbert spaces with  inner product $ \left \langle \cdot, \cdot \right \rangle _{U}$, $\left \langle \cdot,\cdot  \right \rangle _{K}$ and   $(V,\left \| \cdot  \right \| _{V})$ be a reflexive Banach space such that
\begin{eqnarray*}
		V\subset  U= U^{*}\subset V^{*},
	\end{eqnarray*}
where $U^{*}$, $V^{*}$ are the dual spaces of $U$, $V$ and $V\subset U$ continuously and densely.  So
we have $U^{*} \subset V ^{*}$ continuously and densely. Let $_{V ^{\ast}} \langle \cdot , \cdot \rangle _{V} $  denote the pairing between $V^{\ast}$ and $V$, which shows that for all $u \in U$, $v \in V$,
\begin{eqnarray*}
		_{V ^{\ast}} \langle u, v\rangle _{V} = \langle u, v\rangle_{ U},
	\end{eqnarray*}
 and  ($V,U, V ^{\ast}$) is called a Gelfand triple.

 The  $\left ( \Omega ,\mathscr{F},\lbrace\mathscr{F}_{t}\rbrace_{t\ge0},\mathbb{P} \right )$ is a certain complete probability space and $\mathcal{B}(Y)$ denotes the $\sigma $-algebra generated by space $Y$. Let $m\vee n:=\max\{m,n\}$ and $m\wedge n:=\min\{m,n\}$. In this paper, denote by $\mathbb{R}^{n}$ the $n$-dimensional Euclidean space and $\left | \cdot  \right | $ the Euclidean norm. $\mathcal{H}   := C([-\tau , 0]; U)$ is regarded as a space of all  continuous functions from $[-\tau , 0]$ into $U$  and has the norm $\left \| \varphi  \right \| _{\mathcal{H}   } =\underset{-\tau\le \theta \le 0}{\sup}\left \|\varphi (\theta ) \right \|_{U}  $ for all $\varphi \in \mathcal{H} $, where $\tau \in (0,+\infty )$ is referred to as the delay. We use $\mathcal{T}$ to represent the probability measure set on $[-\tau  ,0]$, i.e., for any $\pi \in \mathcal{T} $, $\int_{-\tau  }^{0} \pi  (\text{d}\theta )=1$. For any $q \ge 1$  and  the Banach space $(Y, \left \|   \cdot  \right \|_{Y}   )$,  $\mathcal{L} ^{q}(\Omega, Y) $ denotes the Banach space of all $Y$-value random variables:
\begin{eqnarray*}
		\mathcal{L} ^{q}(\Omega, Y)=\left \{ y:\Omega \to Y :\mathbb{E}\left \| y \right \| ^{q}=\int _{\Omega } \left \| y \right \| ^{p}_{Y}\text{d}\mathbb{P}< \infty   \right \},
	\end{eqnarray*}
where
\begin{eqnarray*}
		\left \| y \right \| _{q} =(\int _{\Omega } \left \| y \right \| ^{q}_{Y}\text{d}\mathbb{P})^{\frac{1}{q} }.
	\end{eqnarray*}
Let $\mathcal{P}(\mathcal{H}) $ be the family of all probability measures on $(\mathcal{H} , \mathcal{B}(\mathcal{H} ))$ with the following bounded Lipschitz distance
 \begin{eqnarray*}
		\left \| \mu _{1}-\mu _{1}\right \|_{BL} :=\sup\{\left | \int F\text{d}\mu _{1}- \int F\text{d}\mu _{2} \right |:\left \| F  \right \|_{BL}\le 1  \},
	\end{eqnarray*}
 where $\left \| F \right \|_{BL}:=\left \| F \right \| _{\infty }+Lip(F)$ and $\left \| F \right \| _{\infty }=\sup_{\varphi \in \mathcal{H} }|F(\varphi)|$, $Lip(F)=\sup_{\varphi_{1}\ne \varphi_{2}}\frac{\left | F(\varphi_{1})-F(\varphi_{2}) \right | }{\left \|\varphi_{1}-\varphi_{2}  \right \|_{\mathcal{H} }} $.
 For a measure-valued map $\mu: \mathbb{R}\to \mathcal{P}(\mathcal{H} )$,  let's further define $\mathcal{R} _{q}$ as follows:
\begin{eqnarray*}
		\mathcal{R} _{q} :=\left \{ \mu_{t} \in  \mathcal{P}(\mathcal{H} ):\underset{t\in \mathbb{R} }{\sup}  \underset{\mathcal{H} }{\int} \|\phi \|^{q} \mu_{t} (\text{d}\phi  ) < \infty   \right \}.
	\end{eqnarray*}

Consider a stochastic functional partial differential equation with finite delay:
\begin{equation}\label{r1}
		\begin{cases}
 \text{d}u(t)=(A(t,u)+f(t,u_{t}))\text{d}t+g(t,u_{t})\text{d}W(t), \\
u_{s} =\varphi\in \mathcal{H},
\end{cases}
	\end{equation}
 where $u_{t} =\{u_{t} (\theta )\}=u(t+\theta )(-\tau \le \theta \le 0)$ and $A(\cdot,\cdot):\mathbb{R}\times V\to V^{*}$ is a family of nonlinear monotone and coercive operators. $f:\mathbb{R}\times\mathcal{H}  \to U$ and $g:\mathbb{R}\times\mathcal{H}  \to \mathscr{L}(K,U)$ are two continuous maps where $\mathscr{L}(K,U)$ is the space of all bounded linear operators from $K$ into $U$. Let $W(t)(t \in \mathbb{R} )$ be a two-sided cylindrical $Q$-Wiener process with $Q = I$ on a separable Hilbert space $(K,\langle \cdot , \cdot \rangle_{ K })$ with respect to a complete
filtered probability space $\left ( \Omega ,\mathscr{F},\lbrace\mathscr{F}_{t}\rbrace_{t\ge0},\mathbb{P} \right )$. To show the dependence of the solution $u(t)$ of system \eqref{r1} on  initial data, we also write $u(t)$ as $u(t;s,\varphi)$.

 We know that the solution $u(t)$ of \eqref{r1} is historically relevant, so it is non-Markov. However, it is proved  that the solution map $u_{t}$ has Markov properties in \cite{R5} and \cite{R10}. So we define the transition probability on space $\mathcal{H}$ for the solution map, i.e., the transition probability of the Markov process defined on $\mathcal{H}$ is a function $p:\Delta \times \mathcal{H} \times \mathcal{B}(\mathcal{H} )\to\mathbb{R} ^{+} $ where $\Delta=\{(t,s):t\ge s, t,s \in \mathbb{R} \}$ for $u_{t} \in \mathcal{H}  $ and $\Gamma\in \mathcal{B}(\mathcal{H} )$ with the following properties:
\begin{enumerate}[1)]
		\item $p(t,s,u_{t},\Gamma)=\mathbb{P}(\omega  :u_{t}\in \Gamma|u_{s})$;
        \item $p(t,s,\cdot,\Gamma)$ is $\mathcal{B}(\mathcal{H} )$-measurable for every $t\ge s$ and $\Gamma\in \mathcal{B}(\mathcal{H} )$;
		\item $p(t,s,\phi ,\cdot)$ is a probability measure on $\mathcal{B}(\mathcal{H} )$ for every $t\ge s$ and $\phi  \in \mathcal{H} $;
        \item The Chapman-Kolmogorov equation:
        \begin{eqnarray*}
		p(t,s,\varphi,\Gamma)=\underset{\mathcal{H} }{\int}p(t,\tau ,\phi,\Gamma)p(\tau,s ,\varphi,\text{d}\phi )
	\end{eqnarray*}
holds for any $s \le \tau \le t  $, $\varphi\in \mathcal{H}$ and $\Gamma\in \mathcal{B}(\mathcal{H} )$.
	\end{enumerate}
We further define a map $\hat{p} (t,s): \mathcal{P}(\mathcal{H} )\to  \mathcal{P}(\mathcal{H} )$ for any $\mu \in \mathcal{P}(\mathcal{H} )$ and $\Gamma\in \mathcal{B}(\mathcal{H} )$ by
\begin{equation}\label{r2}
		\hat{p} (t,s)\mu(\Gamma)=\underset{\mathcal{H} }{\int} p(t,s,\phi,\Gamma)\mu (\text{d}\phi).
	\end{equation}

Next we will give the concept of entrance measures:
\\\textbf{Definition 2.1.} We say a measure-valued map $\mu: \mathbb{R}\to \mathcal{P}(\mathcal{H} )$ is an entrance measure of  \eqref{r1} if $\hat{p} (t,r)\mu_{r}=\mu_{t}$ for all $t\ge r, r\in \mathbb{R}$.
\\\textbf{Lemma 2.2.} \emph{Assume $\mu_{1}$ and $\mu_{2}$ are two probability measures on $(\mathcal{H} , \mathcal{B}(\mathcal{H} ))$. For any open set $\Gamma \subset \mathcal{H} $, if $\mu_{1}(\Gamma)\le\mu_{2}(\Gamma)$,  then $\mu_{1}=\mu_{2}$.
\\\textbf{proof}} The specific proof details can be found in Lemma 2.9 of \cite{R1}.
~\\

Let $(\mathcal{Y},d_{1})$ and $(Z,d _{2})$ be two complete metric spaces; $C(\mathbb{R},Z)$ represent the set of all continuous functions, let  $\Phi ^{l}: =\Phi (t+l)$ be the $l$-translation of $\Phi $,  for any $\Phi \in C(\mathbb{R},Z)$. Let $H(\Phi )$ be the hull of $\Phi $, which is the set of all the limits of $\Phi ^{l_{n} }$ in $C(\mathbb{R},Z)$, i.e., for some sequence $l_{n}\subset \mathbb{R}$
\begin{eqnarray*}
		H(\Phi):=\left \{ \Psi \in C(\mathbb{R},Z) :\lim_{n \to \infty} d_{2}(\Phi ^{l_{n}},\Psi )=0   \right \} .
	\end{eqnarray*}
The specific definitions of various Poisson stable functions and the relationship between their functions can be seen in \cite{R69,R70,R71}.

 Let  $BUC(\mathbb{R}\times\mathcal{Y}   ,Z)$ represent the set of functions satisfying the following properties:
\begin{enumerate}[(1)]
		\item $\Phi$ are continuous in $t$ uniformly w.r.t. $y$ on every bounded subset $B\subseteq \mathcal{Y} $;
         \item $\Phi$ are bounded on all bounded subset from $\mathbb{R}\times\mathcal{Y} $.
	    \end{enumerate}
And we assume that $BUC(\mathbb{R}\times\mathcal{Y} ,Z)$ has the following metric:
\begin{eqnarray*}
		\rho _{BUC} (\Phi,\Psi)=\sum_{i=1 }^{\infty} \frac{1}{2^{n} }\frac{\rho _{n}(\Phi,\Psi) }{1+\rho_{n}(\Phi,g)},
	\end{eqnarray*}
where $\rho _{n}(\Phi,\Psi)=\underset{\left | t \right |\le n,y \in B^{n}  }{\sup} d_{2} (\Phi,\Psi) $, and $\left \{ B^{n}  \right \}$ are bounded, $B^{n}\subset B^{n+1} $   and  $\mathcal{Y} = \underset{n\ge 1}{\bigcup } B^{n} $. Thus $(BUC(\mathbb{R}\times \mathcal{Y} ,Z), \rho_{BUC})$ is a complete metric space(see \cite{R72} for details).

$BC(  \mathcal{Y} ,Z)$ represents the set of all continuous and bounded functions on every bounded subsets  $B\subseteq \mathcal{Y}   $  and have the following metric:
\begin{eqnarray*}
		\rho_{BC} (\Phi,\Psi)=\sum_{i=1 }^{\infty} \frac{1}{2^{n} }\frac{d_{n}(\Phi,\Psi) }{1+d_{n}(\Phi,\Psi)},
	\end{eqnarray*}
 where $d_{n}(\Phi,\Psi) :=\underset{y \in B^{n} }{\sup}\rho (\Phi(y ),\Psi(y )) $. Then  $(BC( \mathcal{D} ,Y), d_{BC})$ is a complete metric space.

 In addition, we need to introduce the following symbols for $\Phi \in C(\mathbb{R},Z)$:
\begin{eqnarray*}
		\mathfrak{N} _{\Phi} =\left \{ \left \{ t_{n}  \right \}\subset \mathbb{R} :\Phi^{t_{n} }\to\Phi   \right \} ,\quad  \mathfrak{N} _{\Phi}^{u} =\left \{\left \{ t_{n}  \right \}\subset \mathfrak{N} _{\Phi} :\Phi^{t_{n} } \; \text{converges} \; \text{to} \; \Phi \;  \text{uniformly}\; \text{in}\; t \in  R \right \},
	   \end{eqnarray*}
\begin{eqnarray*}
		\mathfrak{M} _{\Phi} =\left \{ \left \{ t_{n}  \right \}\subset \mathbb{R}:\left \{ \Phi^{t_{n} }  \right \} \; \text{converges}   \right \} ,\quad  \mathfrak{M} _{\Phi}^{u} =\left \{ \left \{ t_{n}  \right \}\subset \mathfrak{M} _{\Phi}:\left \{ \Phi^{t_{n} }  \right \} \; \text{converges} \;  \text{uniformly}\; \text{in}\; t \in  R \right \}.
	   \end{eqnarray*}
Now for   $f\in BUC(\mathbb{R}\times \mathcal{Y}   ,Z)$, define   $f^{*} :\mathbb{R}\to BC(\mathcal{Y}   ,Z)$ by $f^{*}(t):=f(t,\cdot )$.
\\\textbf{Remark 2.3.} \cite{R69} For all $f\in BUC(\mathbb{R}\times \mathcal{Y}   ,Z)$,
\begin{enumerate}[(1)]
		\item $\mathfrak{M}_{f}= \mathfrak{M}_{f^{*} }$;
		\item $\mathfrak{M}_{f}^{u}= \mathfrak{M}_{f^{*} }^{u}$.
	\end{enumerate}
\textbf{Definition 2.4.} A function $\Phi\in C(\mathbb{R},Z)$ is said to be comparable (by character of recurrence) with $\Psi\in C(\mathbb{R},Z)$ if $\mathfrak{N} _{\Psi} \subseteq \mathfrak{N} _{\Phi}$; $\Phi$ is said to be uniformly comparable (by character of recurrence) with $\Psi$ if $\mathfrak{M} _{\Psi} \subseteq \mathfrak{M} _{\Phi}$.
\\\textbf{Theorem 2.5.} ( \cite{R69,R70,R73}) \emph{ Let $\Phi,\Psi\in C(\mathbb{R},Z)$. Then the following statements hold:
\begin{enumerate}[(1)]
		\item Uniformly comparability implies comparability, i.e., $\mathfrak{M} _{\Psi } \subseteq \mathfrak{M} _{\Phi}$ implies $\mathfrak{N} _{\Psi } \subseteq \mathfrak{N} _{\Phi}$;
		\item $\mathfrak{M} _{\Psi }^{u} \subseteq \mathfrak{M} _{\Phi}^{u}$ implies $\mathfrak{N} _{\Psi }^{u} \subseteq \mathfrak{N} _{\Phi}^{u}$;
         \item Let $\Phi\in C(\mathbb{R},Z)$ be comparable by character of recurrence with $\Psi \in C(\mathbb{R},Z)$. If the function $\Psi $ is stationary (respectively, $T$-periodic, Levitan almost periodic, almost recurrent, Poisson stable), then so is $\Phi$;
         \item Let $\Phi\in C(\mathbb{R},Z)$ be uniformly comparable  by character of recurrence with $\Psi \in C(\mathbb{R},Z)$ and $\Psi$ be Lagrange stable. If $\Psi $ is pseudo-periodic (respectively, pseudo-recurrent), then so is $\Phi$.
	    \end{enumerate}}
~\\
\textbf{Definition 2.6.} Let  $\mu: \mathbb{R}\to \mathcal{P}(\mathcal{H} )$ be a measure-valued map. Then $\mu_{t}$ is called compatible (respectively, uniformly compatible) in the sense of Wasserstein
metric with coefficients, if $\mathfrak{N} _{A,f,g} \subseteq \mathfrak{\widetilde{N} } _{\mu_{t}}$ (respectively, $\mathfrak{M} _{A,f,g} \subseteq \mathfrak{\widetilde{M} } _{\mu_{t}}$), where $\mathfrak{\widetilde{N} } _{\mu_{t}}$ (respectively, $\mathfrak{\widetilde{M} }_{\mu_{t}}$) means the set of all sequences $\left \{ t_{n}  \right \}\subset \mathbb{R}$ such that $\left \| \mu _{(\cdot +t_{n})} -\mu _{(\cdot )}\right \| _{BL}\to 0$ as $t\to \infty $ (respectively, $\{\mu_{(\cdot +t_{n})}\}$ converges)  uniformly on any compact interval.
\section{\textup{Existence and uniqueness of entrance measure}}
Throughout this section, we assume that the initial value $\varphi\in \mathcal{H}$ is independent of $\{W(t)\}_{t\ge s}$. To study the existence and uniqueness of entrance measure  for  \eqref{r1}, we need the following conditions:
 ~\\
 \\\textbf{(A1)} (Boundedness) For $A$, $f$ and $g$, there exist constants $M>0$, $\gamma _{1}>0$ and $p\ge 2$ for all $u\in V$, $t\in\mathbb{R}$ such that
 \begin{eqnarray*}
		\left \| A(t,u)\right \|_{V^{*}}   \le \gamma _{1}\left \| u\right \|_{V}^{p-1}+ M,
	\end{eqnarray*}
and
 \begin{eqnarray*}
		\left \| f(t,0 )\right \|_{U}  \vee \left \| g(t,0 ) \right \|_{\mathscr{L}(K,U)} \le M.
	\end{eqnarray*}
\textbf{(A2)} (Coercivity)  There exist constants $\gamma _{2} \in \mathbb{R}$, $\gamma _{3}>0 $ such that for all $u \in V $, $t \in \mathbb{R}$
\begin{eqnarray*}
		_{V ^{\ast}} \langle A(t,u), u\rangle _{V} \le \gamma _{2}\left \| u \right \| _{H}^{2}-\gamma _{3}\left \| u \right \|_{V}^{p}+M.
	\end{eqnarray*}
\textbf{(A3)} (Monotonicity)  There exists constant $\lambda\in\mathbb{R}$ such that for all $u_{1}, u_{2} \in V $, $t \in \mathbb{R}$
\begin{eqnarray*}
		_{V ^{\ast}} \langle A(t,u_{1})-A(t,u_{2}), u_{1}-u_{2}\rangle _{V} \le \lambda\left \| u_{1}-u_{2} \right \| _{U}^{2}.
	\end{eqnarray*}
\textbf{(A4)} (Semicontinuity) For all $u_{1}, u_{2}, u_{3} \in V $ and $t \in \mathbb{R}$, the map
\begin{eqnarray*}
		\theta  \in \mathbb{R} \to_{V ^{\ast}}\langle A(t,u_{1}+\theta u_{2}),u_{3}\rangle _{V}
	\end{eqnarray*}
is continuous.
\\\textbf{(A5)} For $f$ and $g$, there exist constants $\eta _{1}\in \mathbb{R} $ and $\eta _{2},\eta _{3}, L_{0} >0$ and  $\pi \in \mathcal{T} $, such that for all $\varphi , \phi \in \mathcal{H} $
\begin{eqnarray*}
		\left \langle f(t,\varphi)-f(t,\phi),\varphi (0)-\phi (0) \right \rangle_{U}  \le -\eta _{1} \left \| \varphi (0)-\phi (0) \right \|_{U} ^{2} +\eta _{2}\int_{-\tau  }^{0}\left \| \varphi (\theta )-\phi (\theta ) \right \|^{2}_{U}  \pi  (\text{d}\theta ),
	\end{eqnarray*}
\begin{eqnarray*}
		\left \| f(t,\varphi )-f(t,\phi ) \right \|_{U} \le L_{0}\left \| \varphi -\phi  \right \| _{\mathcal{H} },
	\end{eqnarray*}
\begin{eqnarray*}
		\left \| g(t,\varphi )-g(t,\phi ) \right \| _{\mathscr{L}(K,U)}^{2}  \le \eta _{3}\int_{-\tau  }^{0}\left \| \varphi (\theta )-\phi (\theta ) \right \|^{2}_{U}  \pi  (\text{d}\theta ).
	\end{eqnarray*}

~\\
\\\textbf{Remark 3.1.}
\begin{enumerate}[(1)]
		\item  Under assumptions $\textbf{(A1)-(A5)}$,  the every pair $(\widetilde{A},\widetilde{f} ,\widetilde{g} ) \in H(A,f,g)$ also satisfies the same property with the same constants, where $H(A,f,g)=\overline{\left \{ (A^{l},f^{l},g^{l}  ):l\in \mathbb{R}  \right \} } $ is the hull of $(A,f,g)$;
		\item Under assumptions$\textbf{(A1)-(A5)}$,  $A\in BUC(\mathbb{R} \times V  ,V^{*})$, $f\in BUC(\mathbb{R} \times \mathcal{H}  ,U)$ and  $g\in BUC(\mathbb{R} \times\mathcal{H} ,\mathscr{L}(K,U))$ and $H(A,f,g)\subset  BUC(\mathbb{R} \times V  ,V^{*})\times BUC(\mathbb{R} \times \mathcal{H}  ,U)\times BUC(\mathbb{R} \times\mathcal{H} ,\mathscr{L}(K,U))$.
	\end{enumerate}

In the following, we will present the first important conclusion of the paper and  we will omit the index $U$ of $\| \cdot \|_{U}$ and $\langle \cdot$ ,$ \cdot \rangle_{U}$, if it does not cause confusion.

~\\
\\\textbf{Theorem 3.2.} \emph{Consider \eqref{r1}. Under assumptions \textbf{(A1)}$-$\textbf{(A5)}, and $\lambda<\eta _{1}-\eta _{2}-\frac{145\eta _{3}}{2}  $, there exists a unique entrance measure of \eqref{r1} in $\mathcal{R}_{2}$.
~\\
\\\textbf{proof}} Analyzing system \eqref{r1}, we know that under \textbf{(A1)}$-$\textbf{(A5)},  there exists  a unique solution $u(t;s,\varphi)$ with the initial date $u_{s}=\varphi$ and  a unique solution map $u_{t}(s,\varphi)$ by \cite{R111}.  We divide the proof into five steps:
~\\
\\\textbf{step 1:} There exist  constants $L_{1}$ and $L_{2}$ such that
\begin{eqnarray*}
\begin{split}
		\mathbb{E} \left \| u(t;s,\varphi ) \right \| ^{2} \le L_{1}\left \| \varphi  \right \| _{\mathcal{H}}^{2 }+L_{2},
\end{split}
	\end{eqnarray*}
for all $t\ge s$.

 By It$\hat{\text{o}} $ formula, we obtain for any $\eta>0$, $t\ge s$,
\begin{equation}\label{r4}
\begin{split}
		e^{\eta t} \left \| u(t) \right \| ^{2} &=e^{\eta s}\left \| \varphi (0) \right \| ^{2}+2\int_{s}^{t}e^{\eta r}\left \langle x(r),g(r,u_{r})\text{d}W(r) \right \rangle  \\&~~~+\int_{s}^{t}e^{\eta r} [\eta\left \| u(r) \right \| ^{2}+2_{V ^{\ast}} \langle A(r,u(r)),u(r))\rangle _{V}\\&~~~+2\left \langle f(r,u_{r}),u(r) \right \rangle +\left \| g(r,u_{r}) \right \| ^{2}_{\mathscr{L}(K,U)}  ]\text{d}r.
\end{split}
	\end{equation}
By $\textbf{(A1)}$, $\textbf{(A2)}$, $\textbf{(A5)}$ and Young's inquality, for any $\epsilon _{1},\epsilon_{2}\in (0,1)$, $t\ge s$,
\begin{equation}\label{g13}
\begin{split}
		e^{\eta t} \left \| u(t) \right \| ^{2}&\le e^{\eta s}\left \| \varphi (0) \right \| ^{2} +(\frac{1}{2\epsilon_{1} }+\frac{1}{\epsilon_{2}}  )\frac{M^{2}}{\eta}(e^{\eta t}-e^{\eta s})+2\int_{s}^{t}e^{\eta r}u^{T}(r) g(r,u_{r})\text{d}W(r)\\&~~~+[\eta-(2\eta _{1}-\epsilon_{1})+2\lambda ]\int_{s}^{t}e^{\eta r}\left \| u(r) \right \| ^{2}\text{d}r\\&~~~+(2\eta _{2}+\frac{\eta _{3}}{1-\epsilon_{2}})\int_{s}^{t}e^{\eta r}\int_{-\tau  }^{0}\left \| u(r+\theta ) \right\|^{2}  \pi  (\text{d}\theta )\text{d}r.
\end{split}
	\end{equation}
Further, we have, for $t\ge s$,
\begin{equation}\label{22}
\begin{split}
		&\int_{s}^{t}\int_{-\tau  }^{0}e^{\eta r}\mathbb{E}\left \| u(r+\theta ) \right \|^{2}  \pi  (\text{d}\theta )\text{d}r\\&=\int_{-\tau  }^{0}\int_{s}^{t}e^{\eta r}\mathbb{E}\left \| u(r+\theta ) \right \|^{2} \text{d}r \pi  (\text{d}\theta )\\&\le e^{\eta\tau}\int_{s-\tau }^{t} e^{\eta r}\mathbb{E}\left \| u(r) \right \| ^{2}\text{d}r.
\end{split}
	\end{equation}
Taking expectation of both sides for \eqref{g13}, we have, for $t\ge s$,
\begin{eqnarray*}
\begin{split}
		e^{\eta t} \mathbb{E}\left \| u(t) \right \| ^{2}\le & e^{\eta s}\mathbb{E}\left \| \varphi (0) \right \| ^{2} +(\frac{1}{2\epsilon_{1} }+\frac{1}{\epsilon_{2}}  )\frac{M^{2}}{\eta}(e^{\eta t}-e^{\eta s})+(2\eta _{2}+\frac{\eta _{3}}{1-\epsilon_{2}})e^{\eta \tau}\int_{s-\tau}^{s}e^{\eta r}\mathbb{E}\left \| u(r) \right \| ^{2}\text{d}r\\&+[\eta-(2\eta _{1}-\epsilon_{1})+2\lambda +(2\eta _{2}+\frac{\eta _{3}}{1-\epsilon_{2}})e^{\eta \tau}]\int_{s}^{t}e^{\eta r}\mathbb{E}\left \| u(r) \right \| ^{2}\text{d}r.
\end{split}
	\end{eqnarray*}
Let
\begin{eqnarray*}
\begin{split}
		\Theta (\eta )  =\eta-(2\eta _{1}-\epsilon_{1})+2\lambda +(2\eta _{2}+\frac{\eta _{3}}{1-\epsilon_{2}})e^{\eta \tau}.
\end{split}
	\end{eqnarray*}
Since $\lambda <\eta _{1}-\eta _{2}-\frac{145\eta _{3}}{2} $, there exist $\epsilon_{1},\epsilon_{2}\in (0,1)$ such that $2\eta _{1}-\epsilon_{1}-2\lambda>2\eta _{2}+\frac{\eta _{3}}{1-\epsilon_{2}}$. Hence  $\exists \tilde{\eta } >0$ such that $\Theta (\tilde{\eta } )=0$, which implies
\begin{equation}\label{51}
\begin{split}
		 \mathbb{E}\left \| u(t) \right \| ^{2}\le   L_{1}+L_{2}\left \| \varphi  \right \| _{\mathcal{H}}^{2 },
\end{split}
	\end{equation}
where $L_{1}=(\frac{1}{2\epsilon_{1} }+\frac{1}{\epsilon_{2}}   )\frac{M^{2}}{\tilde{\eta }}$ and $L_{2}=1+(2\eta _{2}+\frac{\eta _{3}}{1-\epsilon_{2}})\frac{e^{\tilde{\eta } \tau}}{\tilde{\eta }}$.
~\\
\\\textbf{step 2:} There exist  constants $L_{3}$,  $L_{4}$ such that
\begin{eqnarray*}
\begin{split}
		\mathbb{E} \left \| u(t;s_{1},\varphi )-u(t;s_{2},\varphi ) \right \| ^{2} \le (L_{3}+L_{4}\left \| \varphi  \right \| _{\mathcal{H}}^{2 })e^{-\eta^{*} (t-s_{2})},
\end{split}
	\end{eqnarray*}
for all $-\infty <s_{2}\le s_{1}\le t<\infty $.

For $-\infty <s_{2}\le s_{1}\le t<\infty $, applying the It$\hat{\text{o}} $ formula to $e^{\eta t} \left \| u(t;s_{1},\varphi )-u(t;s_{2},\varphi ) \right \| ^{2} $ yields
\begin{equation}\label{r8}
\begin{split}
		 &e^{\eta t} \mathbb{E}\left \| u(t;s_{1},\varphi )-u(t;s_{2},\varphi ) \right \| ^{2} \\&=e^{\eta s_{2}}\mathbb{E}\left \| u(s_{2};s_{1},\varphi )-\varphi(0) \right \| ^{2}+\eta\int_{s_{2}}^{t}e^{\eta r} \mathbb{E}\left \|u(r;s_{1},\varphi )-u(r;s_{2},\varphi ) \right \| ^{2}\text{d}r\\&~~~+\int_{s_{2}}^{t}e^{\eta r}\mathbb{E}\left \| g(r,u_{r}(s_{1},\varphi ))-g(r,u_{r}(s_{2},\varphi )) \right \| ^{2}_{\mathscr{L}(K,U)}  \text{d}r\\&~~~+2\mathbb{E}\int_{s_{2}}^{t}e^{\eta r} [_{V ^{\ast}} \langle A(r,u(r;s_{1},\varphi ))-A(r,u(r;s_{2},\varphi )),u(r;s_{1},\varphi )-u(r;s_{2},\varphi ))\rangle _{V}\\&~~~+\left \langle f(r,u_{r}(s_{1},\varphi ))-f(r,u_{r}(s_{2},\varphi )),u(r;s_{1},\varphi )-u(r;s_{2},\varphi ) \right \rangle ]\text{d}r.
\end{split}
	\end{equation}
Similar to \eqref{22},
\begin{equation}\label{r7}
\begin{split}
		&\int_{s_{2}}^{t}\int_{-\tau  }^{0}e^{\eta r}\mathbb{E}\left \| u(r+\theta ;s_{1},\varphi )-u(r+\theta;s_{2},\varphi ) \right \|^{2}  \pi  (\text{d}\theta )\text{d}r\\&\le e^{\eta\tau}\int_{s_{2}-\tau }^{t} e^{\eta r}\mathbb{E}\left \| u(r ;s_{1},\varphi )-u(r;s_{2},\varphi ) \right \| ^{2}\text{d}r.
\end{split}
	\end{equation}
By $\textbf{(A3)}$ and $\textbf{(A5)}$, we have
\begin{equation}\label{r9}
\begin{split}
		&e^{\eta t} \mathbb{E}\left \| u(t;s_{1},\varphi )-u(t;s_{2},\varphi ) \right \| ^{2}\\&\le  e^{\eta s_{2}}\mathbb{E}\left \| u(s_{2};s_{1},\varphi )-\varphi(0)  \right \| ^{2} +(2\eta _{2}+\eta _{3})e^{\eta \tau}\int_{s_{2}-\tau}^{s_{2}}e^{\eta r}\mathbb{E}\left \| u(r;s_{1},\varphi )-\varphi(r) \right \| ^{2}\text{d}r\\&~~~+[\eta-2\eta _{1}+2\lambda +(2\eta _{2}+\eta _{3})e^{\eta \tau}]\int_{s_{2}}^{t}e^{\eta r}\mathbb{E}\left \| u(r;s_{1},\varphi )-u(r;s_{2},\varphi ) \right \| ^{2}\text{d}r.
\end{split}
	\end{equation}
Similarly,  there exists $\eta^{*} >0$ such that $\eta^{*} -2\eta _{1}+2\lambda +(2\eta _{2}+\eta _{3})e^{\eta^{*}  \tau}=0$. Hence
\begin{equation}\label{r10}
\begin{split}
		 \mathbb{E}\left \| u(t;s_{1},\varphi )-u(t;s_{2},\varphi ) \right \| ^{2}\le  (L_{3}+L_{4}\left \| \varphi  \right \| _{\mathcal{H}}^{2 })e^{-\eta^{*} (t-s_{2})},
\end{split}
	\end{equation}
where $L_{3}=L_{1}(1+\frac{(2\eta _{2}+\eta _{3})e^{\eta^{*}  \tau } }{\eta^{*}  } )$ and $L_{4}=(1+L_{2})(1+\frac{(2\eta _{2}+\eta _{3})e^{\eta^{*}  \tau } }{\eta^{*}  } )$.
~\\
\\\textbf{step 3:} There exists  $L_{5}>0$  such that
\begin{eqnarray*}
\begin{split}
		\mathbb{E} \left \| u_{t}(s_{1},\varphi )-u_{t}(s_{2},\varphi ) \right \|_{\mathcal{H}} ^{2} \le L_{5}e^{-\eta^{*} (t-s_{2})},
\end{split}
	\end{eqnarray*}
 for all $-\infty <s_{2}\le s_{1}\le t<\infty $, where $L_{5}=[e^{\eta^{*} \tau }+\eta _{3}e^{2\eta^{*} \tau }+(\frac{(2\eta _{2}+\eta _{3})e^{2\eta^{*} \tau }+36e^{\eta^{*} \tau }}{\eta^{*}})](L_{3}+L_{4}\left \| \varphi  \right \| _{\mathcal{H}}^{2 })$.

By the definition of $\left \| \cdot  \right \|_{\mathcal{H}}$,
\begin{equation}\label{r11}
\begin{split}
		 \mathbb{E}\left \| u_{t}(s_{1},\varphi )-u_{t}(s_{2},\varphi ) \right \|_{\mathcal{H} } ^{2} = \mathbb{E}(\underset{\sigma \in [t-\tau ,t]}{\sup} \left \| u(\sigma  ;s_{1},\varphi )-u(\sigma ;s_{2},\varphi ) \right \| ^{2}).
\end{split}
	\end{equation}
 By $\textbf{(A3)}$, $\textbf{(A5)}$ and \eqref{r7},  we apply It$\hat{\text{o}} $ formula to yield, for $t-\tau\le\sigma$
\begin{equation}\label{r12}
\begin{split}
		 &\left \| u(\sigma ;s_{1},\varphi )-u(\sigma ;s_{2},\varphi ) \right \| ^{2}\\&\le \left \| u(t-\tau;s_{1},\varphi )-u(t-\tau;s_{2},\varphi ) \right \| ^{2}-2(\eta _{1}-\lambda )\int_{t-\tau}^{\sigma }\left\|u(r;s_{1},\varphi )-u(r;s_{2},\varphi )\right\|^{2}\text{d}r\\&~~~+(2\eta _{2}+\eta _{3})\int_{t-2\tau}^{\sigma }\left\|u(r;s_{1},\varphi )-u(r;s_{2},\varphi )\right\|^{2}\text{d}r\\&~~~+2\int_{t-\tau}^{\sigma }\left \langle u(r;s_{1},\varphi )-u(r;s_{2},\varphi ),[g(r,u_{r}(s_{1},\varphi ))-g(r,u_{r}(s_{2},\varphi ))]\text{d}W(r) \right \rangle.
\end{split}
	\end{equation}
By Burkholder–Davis–Gundy inequality, Young's inequality and \eqref{r10}, we have
\begin{equation}\label{r13}
\begin{split}
		 &\mathbb{E}(\underset{\sigma\in [t-\tau ,t]}{\sup}\int_{t-\tau}^{\sigma}\left \langle u(r;s_{1},\varphi )-u(r;s_{2},\varphi ),[g(r,u_{r}(s_{1},\varphi ))-g(r,u_{r}(s_{2},\varphi ))]\text{d}W(r) \right \rangle)\\&\le  6\mathbb{E}(\underset{r\in [t-\tau ,t]}{\sup}\|g(r,u_{r}(s_{1},\varphi ))-g(r,u_{r}(s_{2},\varphi ))\|^{2})^{\frac{1}{2} } (\int_{t-\tau}^{t}\left \|u(r;s_{1},\varphi )-u(r;s_{2},\varphi )\right \|^{2}\text{d}r)^{\frac{1}{2} }\\&\le \frac{\eta _{3} }{2}\mathbb{E} (\underset{r\in [t-\tau ,t]}{\sup}\int_{-\tau }^{0} \left \| u(r+\theta ;s_{1},\varphi )-u(r+\theta;s_{2},\varphi ) \right \|^{2}\pi  (\text{d}\theta ))\\&~~~+18\mathbb{E} (\int_{t-\tau}^{t}\left \|u(r;s_{1},\varphi )-u(r;s_{2},\varphi )\right \|^{2}\text{d}r)\\&\le\frac{\eta _{3} e^{2\eta^{*} \tau }}{2}(L_{3}+L_{4}\left \| \varphi  \right \| _{\mathcal{H}}^{2 })e^{-\eta^{*} (t -s_{2})} +18\mathbb{E} (\int_{t-\tau}^{t}\left \|x(r;s_{1},\varphi )-x(r;s_{2},\varphi )\right \|^{2}\text{d}r).
\end{split}
	\end{equation}
Substituting \eqref{r13} into \eqref{r12}, we have by $\lambda>\eta _{1}-\eta _{2}-\frac{\eta _{3}}{2} $
\begin{eqnarray*}
\begin{split}
		 &\mathbb{E} (\underset{\sigma\in [t-\tau ,t]}{\sup}\left \| u(\sigma;s_{1},\varphi )-u(\sigma;s_{2},\varphi ) \right \| ^{2})\\&\le \mathbb{E} (\underset{\sigma\in [t-\tau ,t]}{\sup}\{\left \| u(\sigma;s_{1},\varphi )-u(\sigma;s_{2},\varphi ) \right \| ^{2} +[2(\eta _{1}-\lambda)-(2\eta _{2}+\eta _{3})]\int_{t-\tau}^{\sigma}\left \|u(r;s_{1},\varphi )-u(r;s_{2},\varphi ) \right \|^{2}\text{d}r\})\\&\le \mathbb{E}\left \| u(t-\tau;s_{1},\varphi )-u(t-\tau;s_{2},\varphi ) \right \| ^{2}+(2\eta _{2}+\eta _{3})\int_{t-2\tau}^{t-\tau}\mathbb{E}\left \|u(r;s_{1},\varphi )-u(r;s_{2},\varphi ) \right \|^{2}\text{d}r\\&~~~+\eta _{3} e^{2\eta^{*} \tau }(L_{3}+L_{4}\left \| \varphi  \right \| _{\mathcal{H}}^{2 })e^{-\eta^{*} (t -s_{2})} +36\mathbb{E} (\int_{t-\tau}^{t}\left \|u(r;s_{1},\varphi )-u(r;s_{2},\varphi ) \right \|^{2}\text{d}r).
\end{split}
	\end{eqnarray*}
By \eqref{r11} and $\textbf{step 2}$, we obtain that
\begin{equation}\label{g9}
\begin{split}
		 &\mathbb{E}\left \| u_{t}(s_{1},\varphi )-u_{t}(s_{2},\varphi ) \right \|_{\mathcal{H} } ^{2} \\&\le [e^{\eta^{*} \tau }+\eta _{3}e^{2\eta^{*} \tau }+(\frac{(2\eta _{2}+\eta _{3})e^{2\eta^{*} \tau }+36e^{\eta^{*} \tau }}{\eta^{*}})](L_{3}+L_{4}\left \| \varphi  \right \| _{\mathcal{H}}^{2 })e^{-\eta^{*} (t -s_{2})}.
\end{split}
	\end{equation}

~\\
\\\textbf{step 4:} If $\lambda<\eta _{1}-\eta _{2}-\frac{145\eta _{3}}{2}  $, then there  exists a unique  map $\mathcal{U} _{t}=\{\mathcal{U} _{t}(\theta)\}:=\{\mathcal{U}(t+\theta):\theta\in[-\tau  ,0]\}\in \mathcal{L} ^{2}(\Omega, \mathcal{H} )$, where $u_{t}(s,\varphi )\overset{\mathcal{L}  ^{2}}{\rightarrow} \mathcal{U}_{t} $ as $s\to -\infty $. From the proof in $\textbf{step 3}$, we obtain that  there  exists $\mathcal{U}_{t}(\varphi )$ such that
\begin{equation}\label{77}
\begin{split}
		\underset{s \to -\infty}{\lim }  \mathbb{E} \left \| u_{t}(s,\varphi )- \mathcal{U} _{t}(\varphi )\right \|_{\mathcal{H} } ^{2}=0,
\end{split}
	\end{equation}
and define
\begin{eqnarray*}
\begin{split}
		\mathcal{U}  _{t}(\varphi ):=\mathcal{L} ^{2}-\underset{s \to -\infty}{\lim }  u_{t}(s,\varphi ).
\end{split}
	\end{eqnarray*}
Similar to  \eqref{r8}-\eqref{r9}, we obtain that for any $\varphi,\phi \in \mathcal{H}$
\begin{equation}\label{r17}
\begin{split}
		 \mathbb{E}\left \| u(t;s,\varphi )-u(t;s,\phi ) \right \| ^{2}&\le (1+\frac{(2\eta _{2}+\eta _{3})e^{\eta^{*} \tau}}{\eta} )\left \| \varphi- \phi \right \| _{\mathcal{H}   }^{2}e^{-\eta^{*} (t-s)}\\&=L_{6}\left \| \varphi- \phi \right \| _{\mathcal{H}   }^{2}e^{-\eta^{*} (t-s)}.
\end{split}
	\end{equation}
Similar to  \eqref{r11}-\eqref{r13} and \eqref{r17}, we deduce
\begin{equation}\label{1}
\begin{split}
		 \mathbb{E}\left \| u_{t}(s,\varphi )-u_{t}(s,\phi ) \right \|_{\mathcal{H}   } ^{2} &\le L_{6}[e^{\eta^{*} \tau }+\eta _{3}e^{2\eta^{*} \tau }+\frac{(2\eta _{2}+\eta _{3})e^{2\tau }+36e^{\tau }}{\eta^{*} } ]\left \| \varphi- \phi \right \| _{\mathcal{H}   }^{2}e^{- \eta^{*}(t -s)}\\&= L_{7}\left \| \varphi- \phi \right \| _{\mathcal{H}   }^{2}e^{- \eta^{*}(t -s)},
\end{split}
	\end{equation}
which implies  that  $\mathcal{U} _{t}(\varphi )$ is independent of the initial date $\varphi \in\mathcal{H}$, and it is denoted by $\mathcal{U} _{t}$.

In addition, according to the definition of the norm $\left \| \cdot  \right \|_{\mathcal{H} }$, we can also get, for any $t\in \mathbb{R}$,
\begin{equation}\label{19}
\begin{split}
		 \mathbb{E} \left \| u_{t} \right \| _{\mathcal{H} }^{2}&=\mathbb{E}(\underset{\theta \in [-\tau ,0]}{\sup} \left \| u(t+\theta ) \right \| )^{2}\\&\le\mathbb{E}(\underset{ u\in [-\tau+s ,t]}{\sup} \left \| u(u) \right \| )^{2}\\&\le \left \| \varphi  \right \| _{\mathcal{H} } ^{2}+\mathbb{E}(\underset{ \sigma \in (s ,t]}{\sup} \left \| u(\sigma ) \right \|^{2} ).
\end{split}
	\end{equation}
Then, applying the It$\hat{\text{o}} $ formula  yields
\begin{equation}\label{20}
\begin{split}
		\mathbb{E}(\underset{ \sigma \in (s ,t]}{\sup} \left \| u(\sigma ) \right \| ^{2})&=\mathbb{E}\left \| \varphi (0) \right \| ^{2}+2\mathbb{E}(\underset{ \sigma \in (s ,t]}{\sup}\int_{s}^{\sigma }\left \langle u(r),g(r,u_{r})\text{d}W(r) \right \rangle )\\&~~~+\mathbb{E}(\underset{ \sigma \in (s ,t]}{\sup}\int_{s}^{\sigma } [2_{V ^{\ast}} \langle A(r,u(r)),u(r))\rangle _{V}\\&~~~+2\left \langle f(r,u_{r}),u(r) \right \rangle +\left \| g(r,u_{r}) \right \| ^{2}_{\mathscr{L}(K,U)}    ]\text{d}r).
\end{split}
	\end{equation}
By $\textbf{(A1)}$, $\textbf{(A2)}$, $\textbf{(A5)}$ and Young's inquality, applying \eqref{22} leads to
\begin{equation}\label{23}
\begin{split}
		&\mathbb{E}(\underset{ \sigma \in (s ,t]}{\sup}\int_{s}^{\sigma } [2_{V ^{\ast}} \langle A(r,x(r)),x(r))\rangle _{V}+2\left \langle f(r,u_{r}),u(r) \right \rangle +\left \| g(r,u_{r}) \right \| ^{2}_{\mathscr{L}(K,U)}   ]\text{d}r)\\&\le\mathbb{E}[\underset{ \sigma \in (s ,t]}{\sup}\{[-(2\eta _{1}-\epsilon_{1})+2\lambda ]\int_{s}^{\sigma }\left \| u(r) \right \| ^{2}\text{d}r+(\frac{1}{2\epsilon_{1} }+\frac{1}{\epsilon_{2}}   )M^{2}( \sigma -s)\\&~~~+(2\eta _{2}+\frac{\eta _{3}}{1-\epsilon_{2}})\int_{s}^{\sigma }\int_{-\tau  }^{0}\left \| u(r+\theta ) \right \|^{2}  \pi  (\text{d}\theta )\text{d}r\}]\\&\le(\frac{1}{2\epsilon_{1} }+\frac{1}{\epsilon_{2}}   )M^{2}(t-s)+(2\eta _{2}+\frac{\eta _{3}}{1-\epsilon_{2}})\int_{s-\tau}^{s}\mathbb{E}\left \| u(r) \right \| ^{2}\text{d}r\\&~~~+[-(2\eta _{1}-\epsilon_{1})+2\lambda +(2\eta _{2}+\frac{\eta _{3}}{1-\epsilon _{2}})]\int_{s}^{t}\mathbb{E}\left \| u(r) \right \| ^{2}\text{d}r\\&=G_{1}( t- s)+G_{2}\int_{s}^{t}\mathbb{E}\left \| u(r) \right \| ^{2}\text{d}r+G_{3},
\end{split}
	\end{equation}
where
\begin{eqnarray*}
\begin{split}
		 G_{1}=(\frac{1}{2\epsilon_{1} }+\frac{1}{\epsilon_{2}}   )M^{2},\quad G_{2}=-2(\eta _{1}-\epsilon_{1})+2\lambda +(2\eta _{2}+\frac{\eta _{3}}{1-\epsilon_{2}}),
\end{split}
	\end{eqnarray*}
\begin{eqnarray*}
\begin{split}
		G_{3}=(2\eta _{2}+\frac{\eta _{3}}{1-\epsilon_{2}})\tau\left \| \varphi  \right \| _{\mathcal{H} }^{2}.
\end{split}
	\end{eqnarray*}
By $\textbf{(A5)}$, Young's inequality, \eqref{22} and applying the Burkholder–Davis–Gundy inequality, we have
\begin{equation}\label{25}
\begin{split}
		&2\mathbb{E}(\underset{ \sigma \in (s ,t]}{\sup}\int_{s}^{\sigma }\left \langle u(r),g(r,u_{r})\text{d}W(r) \right \rangle  )\\&\le12\mathbb{E}(\int_{s}^{t}\left \|u(r)\right \|^{2} \left \|g(r,u_{r})\right \|_{\mathscr{L}(K,U)}^{2} \text{d}r )^{\frac{1}{2} }\\&\le \frac{1}{2} \mathbb{E}(\underset{ \sigma \in (s ,t]}{\sup}\left \|u(\sigma ) \right \|^{2})+72\mathbb{E}[\int_{s}^{t} (2\eta _{3}\int_{-\tau  }^{0}\left \| u(r+\theta ) \right \|^{2}  \pi  (\text{d}\theta )  +2M^{2})\text{d}r]\\&\le\frac{1}{2} \mathbb{E}(\underset{ \sigma \in (s ,t]}{\sup}(\left \|u(\sigma ) \right \|^{2})+G_{4}(t-s)+G_{5}\int_{s}^{t}\mathbb{E}\left \| u(r) \right \| ^{2}\text{d}r+G_{6},
\end{split}
	\end{equation}
where $G_{4}=144M^{2}$, $G_{5}=144\eta _{3}$, $G_{6}=144\eta _{3}\left \| \varphi  \right \|_{\mathcal{H} }^{2}\tau   $.

By  \eqref{19}-\eqref{25}, we have
\begin{equation}\label{29}
\begin{split}
		 \mathbb{E} \left \| u_{t} \right \| _{\mathcal{H} }^{2}&\le3 \left \| \varphi  \right \| _{\mathcal{H} } ^{2}+2(G_{3}+G_{6})+2(G_{1}+G_{4})(t-s)\\&~~~+2(G_{2}+G_{5})\int_{s}^{t}\mathbb{E}\left \| u_{r}\right \|_{\mathcal{H} } ^{2}\text{d}r.
\end{split}
	\end{equation}
Let
\begin{eqnarray*}
\begin{split}
		v (t)= 3 \left \| \varphi  \right \| _{\mathcal{H} } ^{2}+2(G_{3}+G_{6})+2(G_{1}+G_{4})(t-s)+2(G_{2}+G_{5})\int_{s}^{t}v(r)\text{d}r,
\end{split}
	\end{eqnarray*}
which implies that $v(t)$ satisfies the following equation
\begin{eqnarray*}
\begin{split}
		\dot{v }(t) =2(G_{1}+G_{4})+2(G_{2}+G_{5})v(t),
\end{split}
	\end{eqnarray*}
with initial condition $v (s)=3 \left \| \varphi  \right \| _{\mathcal{H} } ^{2}+2(G_{3}+G_{6})$. Solving this equation for $v (t)$, we get
\begin{eqnarray*}
\begin{split}
		v (t)&=[3 \left \| \varphi  \right \| _{\mathcal{H} } ^{2}+2(G_{3}+G_{6})] e^{2(G_{2}+G_{5})(t-s)}-\frac{G_{1}+G_{4}}{G_{2}+G_{5}} (1-e^{2(G_{2}+G_{5})(t-s)}).
\end{split}
	\end{eqnarray*}
We obtain $ \mathbb{E} \left \| u_{t} \right \| _{\mathcal{H} }^{2}\le v (t)$ by comparison principle, and  it follows from the definition of $ \mathbb{E} \left \| u_{t} \right \| _{\mathcal{H} }^{2}$ that for $t\ge s$
\begin{equation}\label{26}
\begin{split}
		 \mathbb{E} \left \| u_{t} \right \| _{\mathcal{H} }^{2}\le[3 \left \| \varphi  \right \| _{\mathcal{H} } ^{2}+2(G_{3}+G_{6})] e^{2(G_{2}+G_{5})(t-s)}-\frac{G_{1}+G_{4}}{G_{2}+G_{5}} (1-e^{2(G_{2}+G_{5})(t-s)}).
\end{split}
	\end{equation}
By $\lambda<\eta _{1}-\eta _{2}-\frac{145\eta _{3}}{2}  $, $\exists \epsilon_{1},\epsilon_{2}$ such that $G_{2}+G_{5}<0$, hence
\begin{equation}\label{g11}
\begin{split}
		\mathbb{E} \left \|\mathcal{U} _{t}  \right \|_{\mathcal{H} }^{2}\le \underset{s \to -\infty}{\lim \sup} \mathbb{E} \left \|u _{t} \right \|_{\mathcal{H} }^{2}=-\frac{G_{1}+G_{4}}{G_{2}+G_{5}}<\infty ,
\end{split}
	\end{equation}
i.e., $\mathcal{U} _{t}\in \mathcal{L} ^{2}(\Omega, \mathcal{H} )$.

For the proof of uniqueness, if there are two $\mathcal{U} _{t}^{1},\mathcal{U} _{t}^{2}\in \mathcal{L} ^{2}(\Omega, \mathcal{H} )$ of \eqref{r1}, by \eqref{1}, we have for any $t\in \mathbb{R}$,
\begin{eqnarray*}
\begin{split}
		\mathbb{E} \left \|\mathcal{U} _{t}^{1}- \mathcal{U} _{t}^{2} \right \|_{\mathcal{H} }^{2}\le L_{7}\left \| \mathcal{U} _{s}^{1}- \mathcal{U} _{s}^{2} \right \| _{\mathcal{H}   }^{2}e^{-\eta^{*} (t -s)}\le L_{7}e^{-\eta^{*} (t -s)}\underset{\tau \in \mathbb{R}}\sup(\mathbb{E} \left \|\mathcal{U} _{\tau }^{1}  \right \|^{2}+\mathbb{E} \left \|\mathcal{U} _{\tau }^{2}  \right \|^{2}).
\end{split}
	\end{eqnarray*}
Let $s\to -\infty $, then $\mathbb{E} \left \|\mathcal{U} _{t}^{1}- \mathcal{U}^{2} _{t}\right \|_{\mathcal{H} }^{2}\to 0$, i.e., $\mathcal{U}_{t}^{1}\overset{\mathcal{L}  ^{2}}=\mathcal{U} _{t}^{2}$.
~\\
\\\textbf{step 5:} There exists a unique entrance measure of \eqref{r1} in $\mathcal{R}_{2}$.

Existence: By \textbf{step 4}, we know that there exists a unique   $\mathcal{U} _{t}\in \mathcal{L} ^{2}(\Omega, \mathcal{H} )$ of \eqref{r1}. Let  $\mu_{t}(\Gamma )=\mathbb{P} (\omega :\mathcal{U} _{t}\in\Gamma )$ and the transition probability $p(t,s,\varphi,\Gamma)=\mathbb{P}(\omega:u_{t}\in\Gamma\mid u_{s})$ for any $\Gamma \in  \mathcal{B}(\mathcal{H} )$, i.e., $\lim_{s \to -\infty}  p(t,s,\varphi,\Gamma)=\mu_{t}(\Gamma )$.

Then for any $\Gamma \in  \mathcal{B}(\mathcal{H} )$ and  $t\ge r, r\in \mathbb{R}$, we have from the
Chapman–Kolmogorov equation of the transition probability that
\begin{equation}\label{2}
\begin{split}
		\hat{p} (t,r)\mu_{r}(\Gamma )&=\underset{\mathcal{H} }{\int} p(t,r,\phi ,\Gamma )\mu _{r}(\text{d}\phi )\\&=\underset{\mathcal{H} }{\int} p(t,r,\phi,\Gamma )\mathbb{P}(\omega :\mathcal{U} _{r}\in \text{d}\phi)\\&=\underset{s \to -\infty}{\lim }\underset{\mathcal{H} }{\int} p(t,r,\phi,\Gamma )p(r,s,\varphi,\text{d}\phi)\\&=\underset{s \to -\infty}{\lim } p(t,s,\varphi,\Gamma )\\&=\mu_{t}(\Gamma ),
\end{split}
	\end{equation}
which implies  $\mu_{t}$ is an entrance measure of \eqref{r1} by \textbf{Definition 2.1}.
~\\

Uniqueness: In the following, we will prove the uniqueness of the entrance measure $\mu_{t}$. Suppose that there exist two entrance measures $\mu,\varrho $ for system \eqref{r1}. By Lemma 2.2, it suffices to prove that for any open set $\Gamma \subset \mathcal{H} $, such that $\mu_{t}(\Gamma)\le\varrho  _{t} (\Gamma)$. For any $t>s$, we obtain
\begin{eqnarray*}
\begin{split}
		\mu_{t}(\Gamma)-\varrho   _{t}(\Gamma)&=\mu_{t}(\Gamma)-\underset{\mathcal{H} }{\int} p(t,s,\varphi ,\Gamma)\varrho   _{s}(\text{d}\varphi )\\&=\underset{\mathcal{H} }{\int} \mu_{t}(\Gamma)-p(t,s,\varphi,\Gamma)\varrho   _{s}(\text{d}\varphi)\\&=\underset{\mathcal{H} }{\int} \mathbb{P}(\omega :\mathcal{U} _{t}\in \Gamma )-\mathbb{P}(\omega :u_{t}(s,\varphi)\in \Gamma)\varrho   _{s}(\text{d}\varphi).
\end{split}
	\end{eqnarray*}
Let $\Gamma_{\rho }:=\{\phi\in \mathcal{H}:\text{dist}(\phi,\Gamma^{c})>\rho \}$, where $\text{dist}(\phi,\Gamma^{c})=\underset{\varphi \in \Gamma ^{c}}{\inf} \left \| \varphi -\phi  \right \|_{\mathcal{H} } $ and it's not hard to get that $\Gamma _{\rho }\subset \Gamma $ and $\Gamma _{\rho }\to \Gamma$ as $\rho \to 0$. Then
\begin{eqnarray*}
\begin{split}
		\mathbb{P}(\omega :\mathcal{U} _{t}\in \Gamma_{\rho })-\mathbb{P}(\omega :\left \| u_{t}(s,\varphi)-\mathcal{U} _{t} \right \| _{\mathcal{H} }\ge\rho )&\le \mathbb{P}(\omega :\mathcal{U} _{t}\in \Gamma_{\rho },\left \| u_{t}(s,\varphi)-\mathcal{U} _{t} \right \| _{\mathcal{H} }<\rho)\\&\le\mathbb{P}(\omega :[u_{t}(s,\varphi)-\mathcal{U} _{t}+\mathcal{U}_{t}]\in \Gamma)\\&= \mathbb{P}(\omega :u_{t}(s,\varphi)\in \Gamma).
\end{split}
	\end{eqnarray*}
By Chebyshev inequality and \eqref{1}, we have
\begin{eqnarray*}
\begin{split}
		 \mathbb{P}(\omega :\mathcal{U} _{t}\in \Gamma )-\mathbb{P}(\omega :u_{t}(s,\varphi)\in \Gamma)&\le \mathbb{P}(\omega :\mathcal{U} _{t}\in \Gamma\setminus \Gamma_{\rho } )+\mathbb{P}(\omega :\left \| u_{t}(s,\varphi)-\mathcal{U}_{t} \right \| _{\mathcal{H} }\ge\rho )\\&\le \mathbb{P}(\omega :\mathcal{U} _{t}\in \Gamma\setminus \Gamma_{\rho } )+\rho^{-2}\mathbb{E}(\left \| u_{t}(s,\varphi)-\mathcal{U} _{t} \right \| _{\mathcal{H} }^{2} )\\&\le \mathbb{P}(\omega :\mathcal{U} _{t}\in \Gamma\setminus \Gamma_{\rho } )+\rho^{-2}L_{7}\left \| \varphi-\mathcal{U} _{s} \right \| _{\mathcal{H}   }^{2}e^{-\eta^{*} (t -s)}\\&\le \mathbb{P}(\omega :\mathcal{U} _{t}\in \Gamma\setminus \Gamma_{\rho } )+\rho^{-2}L_{7}[L_{1}+(2+L_{2})\left \| \varphi \right \| _{\mathcal{H}   }^{2}]e^{-\eta^{*} (t -s)}.
\end{split}
	\end{eqnarray*}
Hence
\begin{eqnarray*}
\begin{split}
		\mu_{t}(\Gamma)-\varrho   _{t}(\Gamma)&=\underset{\mathcal{H} }{\int} \mathbb{P}(\omega :\mathcal{U} _{t}\in \Gamma )-\mathbb{P}(\omega :u_{t}(s,\varphi)\in \Gamma)\varrho   _{s}(\text{d}\varphi)\\&\le \mathbb{P}(\omega :\mathcal{U} _{t}\in \Gamma\setminus \Gamma_{\rho } )+\rho^{-2}L_{7}e^{-\eta^{*} (t -s)}\underset{\mathcal{H} }{\int} [L_{1}+(2+L_{2})\left \| \varphi \right \| _{\mathcal{H}   }^{2}]\varrho _{s}(\text{d}\varphi).
\end{split}
	\end{eqnarray*}
As $s\to -\infty$,
\begin{eqnarray*}
\begin{split}
		\mu_{t}(\Gamma)-\varrho   _{t}(\Gamma)\le \mathbb{P}(\omega :\mathcal{U} _{t}\in \Gamma\setminus \Gamma_{\rho } )=\mu_{t}(\Gamma\setminus \Gamma_{\rho } ).
\end{split}
	\end{eqnarray*}
On the other hand
\begin{eqnarray*}
\begin{split}
		\underset{\rho \to 0 }{\lim}\mu_{t}( \Gamma\setminus \Gamma_{\rho } )=0.
\end{split}
	\end{eqnarray*}
So letting $\rho \to 0 $, we obtain
\begin{eqnarray*}
\begin{split}
		\mu_{t}(\Gamma)-\varrho   _{t}(\Gamma)\le0,
\end{split}
	\end{eqnarray*}
which implies $\mu=\varrho   $ by Lemma 2.2.

Finally, we obtain by \eqref{g11}
\begin{eqnarray*}
\begin{split}
		\underset{t\in \mathbb{R} }{\sup}  \underset{\mathcal{H} }{\int} \|\phi \|^{2} \mu_{t} (\text{d}\phi  )=\underset{t\in \mathbb{R} }{\sup}\mathbb{E} \left \| \mathcal{U} _{t}\right \|_{\mathcal{H}   }^{2}\le-\frac{G_{1}+G_{4}}{G_{2}+G_{5}}<\infty,
\end{split}
	\end{eqnarray*}
i.e., $\mu_{t}\in \mathcal{R}_{2}$.
\quad $\Box $
\section{\textup{Poisson stability entrance measure of SFPDEs with monotone coefficients}}
Based on the above results, in this section we show that the entrance measure $\mu_{t}$ of \eqref{r1} possesses the same recurrent properties as the coefficients $A$, $f$ and $g$. Let $\Lambda \subset \mathbb{R}^{n}(n \in \mathbb{N}) $ be an open bounded subset and $U=U(\Lambda ):=W^{k,2}(\Lambda )$, $k\in \mathbb{N} $.

 In order to analyze the Poisson stability measure, we need the tightness of  the family of $\{\mu _{t}\}_{t\in\mathbb{R} }$. Therefore, we need the following condition $\textbf{(A6)}$, which is used in many works.
~\\
\\\textbf{(A6)} Assume $(\widetilde{U} ,\| \cdot \|_{\widetilde{U}  })   $ is  a closed subset of $(U,\| \cdot \|_{U })   $   such that  $\widetilde{U}   \subset U$ is compact and $V \subset \widetilde{U}  $ is continuous. Let $A_{n}$ be a sequence of positive definite self-adjoint operators on $U$ and for each $n \geq 1$,
\begin{eqnarray*}
\begin{split}
		\left \langle u,v \right \rangle _{n}:=\left \langle u,A_{n}v \right \rangle _{U},\quad u,v\in U,
\end{split}
	\end{eqnarray*}
is a new inner product on $U$. Furthermore, we suppose that the norms $\| \cdot \| _{n}$ generated
by $\langle \cdot , \cdot \rangle_{ n}$ are all equivalent to $\| \cdot \| _{U}$ and we have
\begin{eqnarray*}
\begin{split}
		\left \| u \right \|_{n}\uparrow \left \| u \right \|_{\widetilde{U}  }\quad \text{as} \quad n\to \infty,
\end{split}
	\end{eqnarray*}
 for all $u \in \widetilde{U} $. Let $U_{n} := (U,\left \langle , \right \rangle_{n}  )$ and assume further that  for each $n \geq 1$, $A_{n} : V \rightarrow V $ is continuous and there exist constants $\tilde{\gamma _{2}}$, $\tilde{\eta _{1}}$, $\tilde{\eta _{2}}$, $\tilde{\eta _{3}}$, $\tilde{L_{0}}$, $\widetilde{M} >0$ such that for all $t \in \mathbb{R} $, $u \in V$ and $\varphi , \phi \in \mathcal{H} $
\begin{eqnarray*}
		_{V ^{\ast}} \langle A(t,u), A_{n}u\rangle _{V} \le -\tilde{\gamma _{2}}\left \| u\right \| _{n}^{2}+\widetilde{M},
	\end{eqnarray*}
\begin{eqnarray*}
		\left \langle f(t,\varphi)-f(t,\phi ),\varphi (0)-\phi(0) \right \rangle_{U_{n}}  \le -\tilde{\eta _{1}} \left \| \varphi (0)-\phi(0) \right \|_{n} ^{2} +\tilde{\eta _{2}}\int_{-\tau  }^{0}\left \| \varphi (\theta )-\phi (\theta ) \right \|^{2}_{n}  \pi  (\text{d}\theta ),
	\end{eqnarray*}
\begin{eqnarray*}
		\left \| f(t,\varphi )-f(t,\phi ) \right \|_{U_{n}} \le \tilde{L_{0}}\left \| \varphi -\phi  \right \| _{\tilde{\mathcal{H} }},
	\end{eqnarray*}
\begin{eqnarray*}
		\left \| g(t,\varphi )-g(t,\phi ) \right \| _{\mathscr{L}(K,U_{n})}^{2}  \le \tilde{\eta _{3}}\int_{-\tau  }^{0}\left \| \varphi (\theta )-\phi (\theta ) \right \|^{2}_{n}  \pi  (\text{d}\theta ).
	\end{eqnarray*}
\textbf{Proposition 1.} \emph{Consider \eqref{r1}. Suppose that conditions of Theorem 3.2 and  $\textbf{(A6)}$ hold, then the $\mathcal{U}  _{t}$  satisfies
\begin{equation}\label{r20}
		\sup_{t\in\mathbb{R} } \mathbb{E}\left \| \mathcal{U} (t)\right \|_{\widetilde{U}  }^{2}<\infty \quad \text{and} \quad  \sup_{t\in\mathbb{R} } \mathbb{E}\left \| \mathcal{U}  _{t} \right \|_{\widetilde{\mathcal{H} }   }^{2}<\infty,
	\end{equation}
where $\mathcal{U} (t):=\mathcal{L} ^{2}-\underset{s \to -\infty}{\lim }  u(t;s,\varphi )$ and $\widetilde{\mathcal{H} } =C([-\tau , 0]; \widetilde{U} )$.
\\\textbf{proof}} Similar to the proof of Proposition 1 in \cite{R2} and \eqref{19}-\eqref{26}.
~\\

We are now ready to show the tightness of the set of $\{\mu _{t}\}_{t\in\mathbb{R} }$ of system \eqref{r1}.
~\\
\\\textbf{Lemma 4.1.}  \emph{Suppose that conditions of Theorem 3.2 and $\textbf{(A6)}$  hold. Then for any $\epsilon > 0$ and $T > 0$, there exists a positive constant $R' _{0}$ independent of $\epsilon\in[0, 1]$ such that the solution $u(t;s,\varphi )$ of \eqref{r1} satisfies
\begin{eqnarray*}
		\mathbb{P} (\omega :\{\sup_{r\in [t,t+T]}\left \| u(r;s,\varphi ) \right \|>R  \})<\epsilon, \quad t\ge s,\quad  R\ge R' _{0}.
	\end{eqnarray*}
\textbf{proof}} The specific proof can be found in \textbf{Appendix \uppercase\expandafter{\romannumeral1}}.
~\\
\\\textbf{Lemma 4.2.}  \emph{Suppose that conditions of Theorem 3.2 and $\textbf{(A6)}$  hold. Then $\mathcal{U}  _{t}$ satisfies that for any $\epsilon>0$ and $\kappa >0$, there exists $\delta =\delta (\epsilon,\kappa )>0$ such that
\begin{eqnarray*}
		\mathbb{P} (\omega :\{\sup_{\theta _{1},\theta _{2}\in [-\tau ,0],\left | \theta _{1}-\theta _{2} \right |<\delta  }\left \| \mathcal{U}  _{t}(\theta _{1})-\mathcal{U}  _{t}(\theta _{2}) \right \|\ge \kappa \}  )\le\epsilon, \quad  \forall t\in \mathbb{R}.
	\end{eqnarray*}
\textbf{proof}} The specific proof can be found in \textbf{Appendix \uppercase\expandafter{\romannumeral2}}.
~\\
\\\textbf{Theorem  4.3.} \emph{Suppose that conditions of Theorem 3.2 and $\textbf{(A6)}$ hold. Then $\{\mu _{t}\}_{t\in\mathbb{R} }$ is tight, thus weakly compact in $\mathcal{P}(\mathcal{H} )$.
~\\
\\\textbf{proof}} We just need to prove that for any $\epsilon > 0$, there exists a compact set $\Gamma \subset \mathcal{H}$ such that for all $t\in \mathbb{R} ^{+}$,
 \begin{eqnarray*}
		\mu _{t}(\Gamma )>1-\epsilon.
	\end{eqnarray*}
By \textbf{Proposition 1}, there exists $\hat{R} >0$ such that for any $t\in \mathbb{R}  $,
\begin{equation}\label{r111}
	\mathbb{E} [\sup_{\theta \in [-\tau ,0]} \left \| \mathcal{U}  (t+\theta ) \right \|_{\widetilde{U}  }^{2} ]\le \hat{R},	
	\end{equation}
hence  for any $R>0$ by Chebyshev’s inequality, we obtain
\begin{eqnarray*}
		\mathbb{P}  (\omega :\{\sup_{\theta \in [-\tau ,0]} \left \| \mathcal{U}  (t+\theta ) \right \|_{\widetilde{U}  }^{2}\ge R\})\le \frac{1}{R^{2}}\mathbb{E}[\sup_{\theta \in [-\tau ,0]} \left \| \mathcal{U}  (t+\theta ) \right \|_{\widetilde{U}  }^{2}]\le \frac{\hat{R}}{R^{2}},
	\end{eqnarray*}
 which implies $\mathbb{P}  (\omega :\{\sup_{\theta \in [-\tau ,0]} \left \| \mathcal{U}  (t+\theta ) \right \|_{\widetilde{U}  }^{2}\ge R\})\to 0$ as $R\to \infty   $, i.e., $\forall \epsilon>0$, $\exists R_{1}>0$ such that for all $t\in \mathbb{R} $,
\begin{equation}\label{r112}
	\mathbb{P}  (\omega :\{\sup_{\theta \in [-\tau ,0]} \left \| \mathcal{U}  (t+\theta ) \right \|_{\widetilde{U}  }^{2}\ge R_{1}\})\le \frac{\epsilon}{3} .
	\end{equation}
 Let $\left \{ \Lambda_{n}  \right \} $ be sequences  of $\Lambda$ such that $\Lambda_{n}\subset \Lambda_{n+1} $ for any $n\in \mathbb{N} $ and $\Lambda=\underset{n\ge 1}{\bigcup } \Lambda_{n}$. Since $\mathbb{E} \left \| \mathcal{U} _{t} \right \|_{\mathcal{H} } ^{2}<\infty $, for all $t\in \mathbb{R} $, we have that for $\forall \epsilon>0, m\in \mathbb{N}$, $\exists N_{m}\in \mathbb{N}$ such that
\begin{eqnarray*}
		\mathbb{E} [\sup_{\theta \in [-\tau ,0]} \left \| \mathcal{U}  _{t}(\theta ,x) \right \|_{W^{k,2}(\Lambda\setminus  \Lambda _{N_{m}}) }^{2}]<\frac{\epsilon}{8^{m}},
	\end{eqnarray*}
which implies that for all $t\in \mathbb{R}$ and $m \in \mathbb{N} $
\begin{equation}\label{r113}
	\mathbb{P}(\omega :\sup_{\theta \in [-\tau ,0]}\left \| \mathcal{U}  _{t}(\theta ,x) \right \|_{W^{k,2}(\Lambda\setminus  \Lambda _{N_{m}}) }^{2}\ge \frac{1}{2^{m}} \}) \le \frac{\epsilon }{4^{m}}.
	\end{equation}
From Lemma 4.2, we have that for $\forall \epsilon>0,m\in \mathbb{N}$, $\exists \delta_{m} =\delta (\epsilon,m )>0$ such that
\begin{equation}\label{r115}
		\mathbb{P} (\omega :\{\sup_{\theta _{1},\theta _{2}\in [-\tau ,0],\left | \theta _{1}-\theta _{2} \right |<\delta_{m}  }\left \| \mathcal{U}  _{t}(\theta _{1})-\mathcal{U}  _{t}(\theta _{2}) \right \|> \frac{1}{2^{m}}  \}  )\le\frac{\epsilon}{4^{m}}  , \quad  \forall t\in \mathbb{R}.
	\end{equation}
Denote by for all $m\in\mathbb{N}$,
\begin{eqnarray*}
		\Gamma _{1}=\{\varphi \in \widetilde{\mathcal{H} } :\left \| \varphi  \right \|_{\widetilde{\mathcal{H} } } <R_{1}\},
	\end{eqnarray*}
\begin{eqnarray*}
		\Gamma _{2}^{m}=\{\varphi \in \mathcal{H}:\sup_{\theta \in [-\tau ,0]} \left \| \varphi(\theta ,x) \right \|_{W^{k,2}(\Lambda\setminus  \Lambda _{N_{m}}) }^{2} \le \frac{1}{2^{m}}\},
	\end{eqnarray*}
\begin{eqnarray*}
		\Gamma _{3}^{m}=\{\varphi \in \mathcal{H}:\sup_{\theta _{1},\theta _{2}\in [-\tau ,0],\left | \theta _{1}-\theta _{2} \right |<\delta_{m}  }\left \|\varphi(\theta _{1})-\varphi(\theta _{2}) \right \|\le \frac{1}{2^{m}}\}.
	\end{eqnarray*}
Letting $\Gamma=\Gamma _{1}\cap (\bigcup_{m=1}^{\infty } \Gamma _{2}^{m})\cap (\bigcup_{m=1}^{\infty } \Gamma _{3}^{m})$ implies that $\Gamma$ is precompact in ${\mathcal{H}}$ by the Ascoli-Arzelà theorem,   then we have $\mu _{t}(\Gamma)>1-\epsilon$.
This completes the proof. \quad $\Box $
~\\
\\\textbf{Proposition 2.} \emph{Suppose that \textbf{(A1)}-\textbf{(A5)} hold. Let $u(t;s,\varphi )$ be the solution of \eqref{r1} with initial value $u_{s}=\varphi$. Then there exist  constants $D_{1}$, $D_{2}$ such that
\begin{equation}\label{g1}
		\mathbb{E} \int_{s}^{s+T} \left \| u(t;s,\varphi ) \right \| ^{p}_{V}\text{d}t+\mathbb{E} \int_{s}^{s+T} \left \| A(u(t;s,\varphi ) )\right \| ^{\frac{p}{p-1} }_{V^{*}}\text{d}t\le D_{1}\mathbb{E}\left \| \varphi  \right \| _{\mathcal{H}}^{2 }+D_{2} ,
	\end{equation}
for any $s\in \mathbb{R} $, $T>0$.
\\\textbf{proof}} Similar to the proof of Lemma 3.5 in \cite{R3}.
~\\
\\\\\textbf{Lemma 4.4.} \emph{Suppose that \textbf{(A1)}-\textbf{(A5)} hold. Assume $u^{n}_{t}(s;\varphi )$ and $u_{t}(s;\varphi )$ are the solution maps  of the following system, respectively,
\begin{equation}\label{g2}
		\begin{cases}
 \text{d}u(t)=(A^{n}(t,u)+f^{n}(t,u_{t}))\text{d}t+g^{n}(t,u_{t})\text{d}W(t),\\
u_{s}=\varphi ^{n},
\end{cases}
	\end{equation}
and
\begin{equation}\label{g3}
		\begin{cases}
 \text{d}u(t)=(A(t,u)+f(t,u_{t}))\text{d}t+g(t,u_{t})\text{d}W(t),\\
u_{s}=\varphi.
\end{cases}
	\end{equation}
For any $t\in \mathbb{R} $, $u\in V$ and $\varphi \in \mathcal{H} $, as
$$\lim_{n \to \infty} \left \| A^{n}(t,u)-A(t,u) \right \|_{V^{*}} =0, \quad \lim_{n \to \infty} \left \| f^{n}(t,\varphi )-f(t,\varphi ) \right \| =0,$$
$$\lim_{n \to \infty} \left \| g^{n}(t,\varphi)-g(t,\varphi) \right \|_{ _{\mathscr{L}(K,U)}} =0,$$
we have
\begin{enumerate}[(\textbf{a})]
		\item If $\lim_{n \to \infty} \mathbb{E} \left \| \varphi ^{n}-\varphi  \right \|^{2}_{\mathcal{H} } =0$, then $\lim_{n \to \infty} \mathbb{E} \left \|u ^{n}_{t}(s,\varphi ^{n})-u_{t}(s,\varphi)  \right \|^{2}_{\mathcal{H} } =0$;
\end{enumerate}
\begin{enumerate}[(\textbf{b})]
        \item If $\lim_{n \to \infty}  \left \| \varphi ^{n}-\varphi  \right \|_{\mathcal{H} } =0$ in probability, then $\lim_{n \to \infty} \left \|u ^{n}_{t}(s,\varphi ^{n})-u_{t}(s,\varphi)  \right \|_{\mathcal{H} } =0$ in probability;
	\end{enumerate}
\begin{enumerate}[(\textbf{c})]
        \item If $\lim_{n \to \infty}  \left \|\nu _{\varphi ^{n}}-\nu _{\varphi }  \right \|_{BL } =0$ in $\mathcal{P}(\mathcal{H})$, then $\lim_{n \to \infty} \left \|\nu _{u_{t}^{n}}-\nu _{u_{t}}  \right \|_{BL } =0$ in $\mathcal{P}(C[s,\infty ),\mathcal{H})$, where
            $$\nu _{\varphi ^{n}}(\Gamma)=\mathbb{P} (\omega :\varphi ^{n}\in\Gamma ),\quad \nu _{\varphi }(\Gamma)=\mathbb{P} (\omega :\varphi \in\Gamma ), $$
            $$\nu _{u_{t}^{n}}(\Gamma)=\mathbb{P} (\omega :u_{t}^{n}(s,\varphi ^{n})\in\Gamma ),\quad \nu _{u_{t}}(\Gamma)=\mathbb{P} (\omega :u_{t}(s,\varphi )\in\Gamma ).$$
	\end{enumerate}
\textbf{proof}} Similar to the proof of Lemma 3.5 in \cite{R3}, the above results can be obtained by \eqref{r11}, \eqref{19}, It$\hat{\text{o}} $ formula, Burkholder–Davis–Gundy inequality, Young's inequality, Proposition 2 and Gronwall's lemma.
~\\
\\\textbf{Theorem 4.5.}	\emph{Suppose that conditions of Theorem 3.2 and $\textbf{(A6)}$ hold. Then the unique entrance measure of \eqref{r1} is uniformly compatible.
\\\textbf{proof}}  We now prove that   $\mathfrak{M} _{A,f,g}\subseteq  \mathfrak{\widehat{M} } _{\mu _{t} }$, where  $\mu_{t}(\Gamma )=\mathbb{P} (\omega :\mathcal{U} _{t}\in\Gamma )$. Let $\{t_{n}\}\in \mathfrak{M} _{A,f,g}$, i.e., there exists $(\widehat{A}, \widehat{f} ,\widehat{g}) \in H(A,f,g)$ such that for any $\iota , I_{1}, I_{2}>0$
\begin{eqnarray*}
		\begin{split}
		\lim_{n\to \infty} \underset{\left | t \right |\le \iota ,\left \| u \right \|_{V}\le I_{1}    }{\sup } \left \| A(t+t_{n},u ) -\widehat{A} (t,u )\right \|_{V^{*}}\to 0 , \quad \lim_{n\to \infty} \underset{\left | t \right | \le \iota,\left \| \varphi \right \|_{\mathcal{H} }\le I_{2}   }{\sup } \left \| f(t+t_{n},\varphi ) -\widehat{f} (t,\varphi )\right \|\to 0,
\end{split}
	\end{eqnarray*}
\begin{eqnarray*}
		\begin{split}\lim_{n\to \infty} \underset{\left | t \right | \le \iota,\left \| \varphi \right \|_{\mathcal{H} }\le l_{2}   }{\sup } \left \| g(t+t_{n},\varphi ) -\widehat{g} (t,\varphi )\right \|_{\mathscr{L}(K,U)}\to 0.
\end{split}
	\end{eqnarray*}
Assume $\mu ^{n}_{t}$ and $\tilde{\mu }_{t} $ are the entrance measure   of the following system, respectively,
\begin{equation}\label{g5}
\text{d}u(t)=(A(t+t_{n},u)+f(t+t_{n},u_{t}))\text{d}t+g(t+t_{n},u_{t})\text{d}W(t),
	\end{equation}
and
\begin{equation}\label{g6}		\text{d}u(t)=(\widehat{A}(t,u)+\widehat{f}(t,u_{t}))\text{d}t+\widehat{g}(t,u_{t})\text{d}W(t).
	\end{equation}
Let $\ell \le-1(\left | \ell  \right | \in \mathbb{N} )$, and for any sequence $\{n:n\in \mathbb{N} \}$, we obtain that  $\{\mu ^{n}_{\ell }\}_{n }$ is  tight by Theorem  4.3. Hence  there exists a subsequence $\{n_{k}\}\subset \{n\}$ such that
\begin{eqnarray*}
		 \mu ^{n_{k}}_{\ell } \overset{W}{\rightarrow}\varrho _{\ell } \quad \text{as}\quad  k\to\infty \quad \text{in}\quad \mathcal{P} (\mathcal{H} ),
	\end{eqnarray*}
where $\varrho _{\ell }\in \mathcal{P} (\mathcal{H} )$. Assume  $u(t,\ell, \varphi _{\ell })$ is the solution  of the following system
\begin{eqnarray*}
		\begin{cases}
 \text{d}u(t)=(\widehat{A}(t,u)+\widehat{f}(t,u_{t}))\text{d}t+\widehat{g}(t,u_{t})\text{d}W(t),\\
u_{\ell}=\varphi _{\ell }\in \mathcal{H} ,
\end{cases}
	\end{eqnarray*}
where $\varphi _{\ell }$ is the random variable with distribution $\varrho _{\ell }$. Let $\nu ^{\ell }_{t}(\Gamma )=\mathbb{P} (\omega :u_{t}(\ell, \varphi _{\ell })\in\Gamma )$. By (c) of Lemma 4.4, we obtain for any $t\in [\ell,\infty )$
\begin{eqnarray*}
		 \lim_{k \to \infty} \left \| \mu ^{n_{k}}_{t }- \nu ^{\ell }_{t}\right \| _{BL}=0\quad \text{in}\quad \mathcal{P} (\mathcal{H} ).
	\end{eqnarray*}
Let $\ell=\ell-1$, and we obtain that $\{\mu ^{n_{k}}_{\ell-1 }\}_{k }$ is  tight by Theorem  4.3. Thus, there exists a subsequence $\{n_{k_{L}}\}:=\{n_{k}\}$ such that $\{\mu ^{n_{k_{L}}}_{\ell-1 }\}_{L }$ is convergent. Going if necessary to a subsequence, we assume the subsequence $\{n_{k_{L}}\}:=\{n_{k}\}$. Repeating the above process, in fact, we can obtain for any $t\in [\ell-1,\infty )$
\begin{eqnarray*}
		 \lim_{k \to \infty} \left \| \mu ^{n_{k}}_{t }- \nu ^{\ell-1 }_{t}\right \| _{BL}=0\quad\text{in}\quad \mathcal{P} (\mathcal{H} ),
	\end{eqnarray*}
which implies that for any $t\in [\ell,\infty )$
\begin{eqnarray*}
		 \nu ^{\ell }_{t}= \nu ^{\ell-1 }_{t}\quad\text{in}\quad \mathcal{P} (\mathcal{H} ).
	\end{eqnarray*}
Note that  $\nu ^{\ell }_{t}$ is independent of $\ell$, which is denoted by $\nu_{t}$. Hence by the standard diagonal argument,  there exists a subsequence $\{n_{k}\}\subset \{n\}$ such that
\begin{eqnarray*}
		 \lim_{k \to \infty} \left \| \mu ^{n_{k}}_{t }- \nu_{t}\right \| _{BL}=0\quad \text{in}\quad \mathcal{P} (\mathcal{H} ).
	\end{eqnarray*}
\\Let $\ell\to -\infty $, then  by  \eqref{77},  there  exists $\mathcal{\widehat{U}}_{t}\in \mathcal{L} ^{2}(\Omega, \mathcal{H} )$ such that
\begin{eqnarray*}
\underset{\ell \to -\infty}{\lim }  \mathbb{E} \left \| u_{t}(\ell,\varphi _{\ell})- \mathcal{\widehat{U}}_{t}\right \|_{\mathcal{H} } ^{2}=0.
	\end{eqnarray*}
 Since $\mathcal{L} ^{2}$-convergence implies convergence in distribution, we have $u_{t}(\ell,\varphi _{\ell})\to \mathcal{\widehat{U}}_{t}$ in distribution uniformly on $\mathbb{R}$, i.e., $\nu_{t}=\widehat{\mu }_{t} $, where $\widehat{\mu }_{t}(\Gamma )=\mathbb{P} (\omega :\mathcal{\widehat{U}}_{t}\in\Gamma )$. Hence we can extract a subsequence which  we still denote by $\mu ^{n_{k}}_{t }$ satisfying
 \begin{eqnarray*}
		 \lim_{k \to \infty} \left \| \mu ^{n_{k}}_{t }- \widehat{\mu }_{t} \right \| _{BL}=\lim_{k \to \infty} \left \| \mu ^{n_{k}}_{t }- \nu _{t}\right \| _{BL}=0\quad \text{in}\quad \mathcal{P} (\mathcal{H} ),
	\end{eqnarray*}
 for any $t\in\mathbb{R} $. The above proof process shows that for  sequence $\{n_{k}\}_{k=1}^{\infty }\subset \mathbb{N} $, there is a subsequence $\{n_{k_{L}}\}_{L=1}^{\infty }\subset\{n_{k}\}_{k=1}^{\infty }$ such that for any $t\in\mathbb{R} $
 \begin{eqnarray*}
		 \lim_{L \to \infty} \left \| \mu ^{n_{k_{L}}}_{t }- \widehat{\mu }_{t} \right \| _{BL}=0\quad \text{in}\quad \mathcal{P} (\mathcal{H} ),
	\end{eqnarray*}
 which implies
 \begin{equation}\label{g7}
		 \lim_{n \to \infty} \left \| \mu ^{n}_{t }- \widehat{\mu }_{t} \right \| _{BL}=0\quad \text{in}\quad \mathcal{P} (\mathcal{H} ).
	\end{equation}

In addition, under \textbf{(A1)}-\textbf{(A5)},  \eqref{g5} admits a unique solution $u^{n}(t;s, \varphi)$ with the initial data $u_{s}=\varphi$ satisfying
\begin{equation}\label{3}
\begin{split}
		u(t;s,\varphi)&=\varphi(0)+\int_{s}^{t}[A(r+t_{n},u(r;s,\varphi))+ f(r+t_{n},u_{r}(s,\varphi ))]\text{d}r\\&~~~+\int_{s}^{t} g(r+t_{n},u_{r}(s,\varphi ))\text{d}W(r),
\end{split}
	\end{equation}
and
\begin{equation}\label{4}
\begin{split}
		&u(t+t_{n};s+t_{n},\varphi)\\&=\varphi(0)+\int_{s+t_{n}}^{t+t_{n}}[A(r,u(r;s+t_{n},\varphi))+ f(r,u_{r}(s+t_{n},\varphi ))]\text{d}r\\&~~~+\int_{s+t_{n}}^{t+t_{n}} g(r,u_{r}(s+t_{n},\varphi ))\text{d}W(r)\\&=\varphi(0)+\int_{s}^{t}[A(r+t_{n},u(r+t_{n};s+t_{n},\varphi))+ f(r+t_{n},u_{r+t_{n}}(s+t_{n},\varphi ))]\text{d}r\\&~~~+\int_{s}^{t} g(r+t_{n},u_{r+t_{n}}(s+t_{n},\varphi ))\text{d}\widehat{W} (r),
\end{split}
	\end{equation}
where $\widehat{W} (t)=W(t+T)-W(t)$ is a  two-sided cylindrical Wiener process with the same distribution as $W(t)$.

Comparing equations \eqref{4} with \eqref{3} and noting that $u^{n}_{t}(s,\varphi)$ and $u_{t+t_{n}}(s+t_{n},\varphi)$ completely depends on $u^{n}(t)$ and $u(t+t_{n})$ and their history, we see by the weak uniqueness that $u^{n}_{t}$ and $u_{t+t_{n}}$ share the same distribution for any $t\ge s$  $(s\in \mathbb{R}) $ . Furthermore, from the conclusion of \textbf{Theorem 3.2}, we can obtain that
\begin{equation}\label{59}
\begin{split}
		\mathcal{U}  _{t+t_{n}}=\mathcal{L} ^{2}-\underset{s \to -\infty}{\lim }  u_{t+t_{n}}(s+t_{n},\varphi )\overset{\mathbb{P} }{=} \mathcal{L} ^{2}-\underset{s \to -\infty}{\lim }  u^{n}_{t}(s,\varphi )=\mathcal{U} ^{n} _{t},
\end{split}
	\end{equation}
i.e., for any $ \Gamma \in \mathcal{B}(\mathcal{H} )$,
\begin{equation}\label{g8}
\begin{split}
		\mu_{t+t_{n}}( \Gamma )=\mathbb{P}(\omega : \mathcal{U} _{t+t_{n}}\in \Gamma )=\mathbb{P}(\omega : \mathcal{U} ^{n}_{t}\in \Gamma )=\mu^{n}_{t}( \Gamma ).
\end{split}
	\end{equation}
Hence, by \eqref{g7} and \eqref{g8}, we have
\begin{eqnarray*}
		 \lim_{n \to \infty} \left \| \mu_ {t+t_{n} }- \widehat{\mu }_{t}\right \| _{BL}=0\quad \text{in}\quad \mathcal{P} (\mathcal{H} ),
	\end{eqnarray*}
for any $t\in\mathbb{R} $, which implies  that  for any $\iota>0$,
\begin{eqnarray*}
		\begin{split}
		\lim_{n \to \infty} \underset{\left | t \right |\le \iota }{\max }   \left \|  \mu_ {t+t_{n} }- \widehat{\mu }_{t} \right \|_{BL} =0,
\end{split}
	\end{eqnarray*}
i.e., the  measure  $ \mu_ {t}$ of \eqref{r1} is uniformly compatible by Definition 2.6.
\quad $\Box $
~\\
\\\textbf{Corollary 4.6.} \emph{Consider \eqref{r1}. Assume  the conditions of Theorem 4.5 hold. By Theorem 2.5, we obtain:
\begin{enumerate}[1)]
        \item  If $A$, $f$ and $g$ are jointly stationary (respectively, $T$-periodic,   almost periodic, Bohr almost automorphic, Birkhoff recurrent, Lagrange stable, Levitan almost periodic, almost recurrent, Poisson stable) in $t\in \mathbb{R} $ uniformly with respect to  $u\in V$ and $\varphi \in\mathcal{H} $ on every bounded subset,  we can obtain that the unique entrance measure $\mu_{t}$ of \eqref{r1} is stationary (respectively, $T$-periodic, almost periodic, Bohr almost automorphic, Birkhoff recurrent, Lagrange stable, Levitan almost periodic, almost recurrent, Poisson stable);
		\item  If $A$, $f$ and $g$ are jointly pseudo-periodic (respectively, pseudo-recurrent) and $A$, $f$ and $g$ are jointly Lagrange stable in $t\in \mathbb{R} $ uniformly with respect to  $u\in V$ and $\varphi \in\mathcal{H} $ on every bounded subset, then the unique entrance measure $\mu_{t}$ of \eqref{r1} is pseudo-periodic (respectively, pseudo-recurrent).
	\end{enumerate}
\textbf{proof}} These statements follow from Definition 2.6, Theorems 2.5, 3.2 and 4.5.

\section{\textup{Exponential mixing }}
Next we will study the autonomous version of SFPDEs:
\begin{equation}\label{rrr1}
		\begin{cases}
 \text{d}u(t)=(A(u)+f(u_{t}))\text{d}t+g(u_{t})\text{d}W(t), \\
u_{s} =\varphi\in \mathcal{H},
\end{cases}
	\end{equation}
 where $A(\cdot,\cdot): V\to V^{*}$ is a family of nonlinear monotone and coercive operators. $f:\mathcal{H}  \to U$ and $g:\mathcal{H}  \to \mathscr{L}(K,U)$ are two continuous maps. Suppose that Hypotheses $\textbf{(A1)-(A6)}$ are satisfied. Therefore, Theorem 3.2 and Theorem 4.3 are still valid for autonomous systems.

Based on the tightness of the set of measures $\{\mu _{t}\}_{t\in\mathbb{R} }$ as established in Section 4, we now investigate  the exponential mixing of \eqref{rrr1}.  With the transition probability $p(t, s, \phi, \Gamma  ) :=\mathbb{P}  (\omega :u_{t}( s, \varphi)\in\Gamma )$, for any $r\ge0$, we associate a mapping $P^{*}_{r}:\mathcal{P}(\mathcal{H} )\to  \mathcal{P}(\mathcal{H} )$ defined by
\begin{equation}\label{r116}
		P^{*}_{r}\mu _{t}(\Gamma )=\hat{p}(r,0)\mu _{t}(\Gamma )=\int _{\mathcal{H} }p(r,0,\phi ,\Gamma )\mu _{t}(\text{d}\phi ).
	\end{equation}
For any $F\in C_{b}(\mathcal{H})$, which is defined as the the set of all bounded continuous functions $F : \mathcal{H}\to \mathbb{R}$ endowed with the norm $\left \| F \right \| _{\infty }=\sup_{\xi \in \mathcal{H} } \left | F(\xi ) \right | $, we define the following semi-group $P_{s,t}$ for $t\ge s$
\begin{equation}\label{r117}
		P_{s,t}F(\xi  )=\int _{\mathcal{H} }F(\phi )p(s,t,\xi   ,\text{d}\phi  ).
	\end{equation}
In particular, the  operator $P_{0,t}$ is written as $P_{t}$.
~\\
\\\textbf{Proposition 3.}\emph{ Suppose that conditions of Theorem  4.5 hold. Then the semi-group  $P_{s,t}(t\ge s)$ is Feller, i.e. for all $F\in C_{b}(\mathcal{H})$, $P_{s,t}F\in C_{b}(\mathcal{H})$.
~\\
\\\textbf{proof}} From the definition of the semi-group $P_{s,t}$, we have $\left \| P_{s,t}F \right \| _{\infty}\le\left \| F \right \| _{\infty }$. The next major task is to prove that $ P_{s,t}F$ is continuous.

We just need to prove that for any sequence $\xi   _{n}\in \mathcal{H} $, $\xi  \in \mathcal{H} $, when $\lim_{n \to \infty } \left \| \xi  _{n}-\xi   \right \| _{\mathcal{H} }=0$, we have $\lim_{n \to \infty }  | P_{s,t}F(\xi  _{n})- P_{s,t}F(\xi  )   |=0$. Since
\begin{equation}\label{r118}
\begin{split}
		P_{s,t}F(\xi   )&=\int _{\mathcal{H} }F(\phi )p(t,s,\xi   ,\text{d}\phi  )\\&
=\int _{\mathcal{H} }F(\phi )\mathbb{P} (\omega :u_{t}(s,\xi  )\in\text{d}\phi  )\\&
=\mathbb{E} F(u_{t}(s,\xi  )),
\end{split}
	\end{equation}
i.e.
\begin{equation}\label{r119}
\begin{split}
		|P_{s,t}F(\xi   _{n} )-P_{s,t}F(\xi   )|=|\mathbb{E} F(u_{t}(s,\xi  _{n}))-\mathbb{E} F(u_{t}(s,\xi  ))|,
\end{split}
	\end{equation}
by \eqref{1} we can get
\begin{eqnarray*}
		\mathbb{E}\left \| u_{t}(s,\xi   )-u_{t}(s,\xi  _{n} ) \right \|_{\mathcal{H}   } ^{2} \le  L_{7}\left \| \xi  - \xi  _{n} \right \| _{\mathcal{H}   }^{2}e^{-\eta^{*} (t-s) },
	\end{eqnarray*}
i.e. $u_{t}(s,\xi  _{n})\overset{L_{2}}{\rightarrow}u_{t}(s,\xi   )$ as $n\to \infty  $. Let
$$W_{R}=\{\omega :\left \|  u_{t}(s,\xi   ,\omega )\right \| _{\mathcal{H} }\le R\},\quad W_{R}^{n}=\{\omega :\left \| u_{t}(s,\xi  _{n} ,\omega )\right \| _{\mathcal{H} }\le R\}.$$
Then by the Chebyshev inequality and \eqref{26} we have
\begin{eqnarray*}
		\lim_{R\to \infty} [\inf_{n\in \mathbb{N} } \mathbb{P}(W_{R}\cap W_{R}^{n} ) ]=1 ,
	\end{eqnarray*}
i.e. $\forall \varepsilon>0$, $\exists \tilde{R} $ such that as $R>\tilde{R}$, for any $n\in \mathbb{N}$
\begin{eqnarray*}
		\mathbb{P}(W_{R}\cap W_{R}^{n} )>1-\frac{\varepsilon}{4\left \| F \right \|_{0}} .
	\end{eqnarray*}
Since  $F\in C_{b}(\mathcal{H})$, we have that $\forall \varepsilon>0$ , $\exists \delta =\delta (\varepsilon,R) $ such that as $\left \| \xi  _{1} -\xi  _{2}  \right \| _{\mathcal{H} }<\delta $, then
\begin{eqnarray*}
		\left | F(\xi  _{1} )-F(\xi  _{2} ) \right | <\varepsilon.
	\end{eqnarray*}
Let $W_{\delta }^{n}=\{\omega :\left \| u_{t}(s,\xi  ,\omega  )-u_{t}(s,\xi  _{n},\omega  ) \right \| <\delta \}$ and we obtain that $\lim_{n\to \infty}  \mathbb{P}(W_{\delta }^{n}) =1$.

For any $\omega \in W_{R}\cap W_{R}^{n}\cap W_{\delta }^{n}$,  we obtain
\begin{eqnarray*}
		|F( u_{t}(s,\xi  ))-F(u_{t}(s,\xi _{n}))|<\frac{\varepsilon}{2}  .
	\end{eqnarray*}
Then
\begin{equation}\label{r119}
\begin{split}
		&\underset{n\to\infty }{\lim \sup } \left | \mathbb{E}F(u_{t}(s,\xi _{n}  ))-\mathbb{E}F(u_{t}(s,\xi  ))   \right |\\&\le\frac{\varepsilon}{2} +2\left \| F \right \|_{\infty } \underset{n\to\infty }{\lim \sup } [2-(\mathbb{P}(W_{R}\cap W_{R }^{n}) +\mathbb{P}(W_{\delta }^{n}) )].
\end{split}
	\end{equation}
Letting $\varepsilon\to 0$ implies $\lim_{n \to \infty }  | P_{s,t}F(\xi _{n})- P_{s,t}F(\xi )   |=0$. Then the semi-group  $P_{s,t}$ is Feller. \quad $\Box $
~\\

For the above entrance measure $\mu_{t}$, set
\begin{equation}\label{r222}
\begin{split}
		\mu ^{L}=\frac{1}{L} \int_{0}^{L} \mu _{t}\text{d}t,
\end{split}
	\end{equation}
and $\mathfrak{L} :=\left \{ \mu ^{L};L\in \mathbb{N}^{+} \right \} $. By Theorem  4.3, we obtain that $\{\mu _{t}\}_{t\in \mathbb{R}}$ is tight on $\mathcal{H}$, which implies that for any $\epsilon > 0$, there exists a precompact set $\Gamma \subset \mathcal{H}$ such that for all $L\in \mathbb{N} ^{+}$, we have
 \begin{eqnarray*}
		\mu ^{T}(\Gamma )=\frac{1}{L} \int_{0}^{L} \mu _{t}(\Gamma)\text{d}t>1-\epsilon,
	\end{eqnarray*}
i.e., $\mathfrak{L}$ is tight, hence $\mathfrak{L}$ is weakly compact in $\mathcal{P}(\mathcal{H} )$. Thus we know that there exists a probability measure $\mu^{*} $ on $\mathcal{P(\mathcal{H}) }  $ such that $\mu ^{L}\to \mu^{*}$ as $L\to \infty $. Before proving the exponential mixing of measure $\mu^{*}$, we still need the following lemma:
~\\
\\\textbf{Lemma 5.1.} \emph{Assume the conditions of Theorem 4.5 hold. Then the measure $\mu^{*}$  satisfies
\begin{eqnarray*}
		P^{*}_{r}\mu^{*}(\Gamma)=\mu^{*}(\Gamma),
	\end{eqnarray*}
for any $r\ge0$ and $\Gamma\in \mathcal{B}(\mathcal{H} )$. In particular, the measure satisfying the above property is unique.
~\\
\\\textbf{proof}} For every $r\ge0$, $t\in \mathbb{R} $ and $\Gamma\in \mathcal{B}(\mathcal{H} )$, by  Proposition 3, Chapman–Kolmogorov equation and Lemma 2.3, we find that
\begin{equation}\label{r33}
\begin{split}
		P_{r}^{*}\mu _{t}(\Gamma )&=\int _{\mathcal{H} }p(r,0,\phi ,\Gamma )\mu _{t}(\text{d}\phi )\\&=\lim_{s \to -\infty}\int _{\mathcal{H} }p(r,0,\phi ,\Gamma )p(t,s,\varphi  ,\text{d}\phi )\\&=\lim_{s \to -\infty}\int _{\mathcal{H} }p(r+t,t,\phi ,\Gamma )p(t,s,\varphi  ,\text{d}\phi ) \\&=\mu _{r+t}(\Gamma ).
\end{split}
	\end{equation}
Therefor, by \eqref{r33} we have
\begin{equation}\label{r36}
\begin{split}
		P_{r}^{*}\mu ^{L}-\mu ^{L}=\frac{1}{L} \int_{0}^{L} \mu_{r+t}\text{d}t-\frac{1}{L} \int_{0}^{L} \mu_{t}\text{d}t=\frac{1}{L} [\int_{L}^{r+L} \mu_{t}\text{d}t+\int_{0}^{r} \mu_{t}\text{d}t],
\end{split}
	\end{equation}
which implies $\lim_{L \to \infty} \left \| P_{r}^{*}\mu ^{L}-\mu^{*}  \right \| _{BL}=0$, i.e. $P_{r}^{*}\mu ^{L} \to\mu^{*} $ weakly in $\mathcal{P}(\mathcal{H})$ as $L \to \infty$.

In addition, for any $F\in C_{b}(\mathcal{H})$, we have $P_{t}F\in C_{b}(\mathcal{H})$ by Proposition 3. Hence by Lebesgue dominated convergence theorem, we have
\begin{equation}\label{r35}
\begin{split}
		\lim_{L \to \infty} \int _{\mathcal{H} }F(\phi )P_{r}^{*}\mu ^{L}(\text{d}\phi )&=\lim_{L \to \infty} \int _{\mathcal{H} }\int _{\mathcal{H} }F(\phi )p(r,0,\varphi ,\text{d}\phi )\mu ^{L}(\text{d}\varphi )\\&=\lim_{L \to \infty} \int _{\mathcal{H} }P_{r}F(\varphi )\mu ^{L}(\text{d}\varphi )\\&=\int _{\mathcal{H} }P_{r}F(\varphi)\mu^{*}(\text{d}\varphi )\\&= \int _{\mathcal{H} }\int _{\mathcal{H} }F(\phi )p(r,0,\varphi ,\text{d}\phi )\mu^{*} (\text{d}\varphi )\\&=\int _{\mathcal{H} }F(\phi)P_{r}^{*}\mu^{*}(\text{d}\phi ),
\end{split}
	\end{equation}
which implies  $P_{r}^{*}\mu ^{L} \to P_{r}^{*}\mu^{*} $ weakly in $\mathcal{P}(\mathcal{H})$.  Therefore  we obtain $P_{r}^{*}\mu ^{*}=\mu^{*} $   by \eqref{r36} and \eqref{r35}.

Uniqueness: Suppose that there exists the measure ${\mu }'\in \mathcal{P} (\mathcal{H} ) $ such that ${\mu }'(\Gamma )=P_{r}{\mu }'(\Gamma )$, then  for any open set $\Gamma \in \mathcal{B} (\mathcal{H} )$, we have
\begin{equation}\label{r37}
\begin{split}
		{\mu }'(\Gamma )=P_{r}{\mu }'(\Gamma )&=\lim_{L \to \infty} \frac{1}{L} \int_{0}^{L} P_{r}{\mu }'(\Gamma )\text{d}t\\&=\lim_{L \to \infty} \frac{1}{L} \int_{0}^{L}\int _{\mathcal{H} }p(t,0,\varphi ,\Gamma ){\mu }'(\text{d}\varphi  )\text{d}t\\&=\lim_{L \to \infty} \frac{1}{L} \int _{\mathcal{H} }\int_{0}^{L}\mathbb{P} (\omega :u_{t}(0,\varphi )\in \Gamma )\text{d}t{\mu }'(\text{d}\varphi  ).
\end{split}
	\end{equation}
By the globally asymptotic  stability of $\mathcal{U}_{t}$, we have
\begin{eqnarray*}
		\mathbb{E} \left \| u_{t}(0,\varphi ) -\mathcal{U} _{t}\right \|^{2}_{\mathcal{H}   } \le L_{7}\left \| \varphi-\mathcal{U} _{0} \right \| _{\mathcal{H}   }^{2}e^{-\eta^{*} t }\le \widetilde{L _{7}} e^{-\eta^{*} t }.
	\end{eqnarray*}
By the Chebyshev inequality, we have for all $\rho >0$
\begin{equation}\label{r38}
\begin{split}
		\mathbb{P}(\omega :u_{t}(0,\varphi)\in \Gamma)&\ge \mathbb{P}(\omega :\mathcal{U} _{t}\in \Gamma_{\rho },\left \| u_{t}(0,\varphi)-\mathcal{U} _{t} \right \| _{\mathcal{H} }<\rho )\\&\ge\mathbb{P}(\omega :\mathcal{U}_{t}\in \Gamma_{\rho })-\mathbb{P}(\omega :\left \| u_{t}(0,\varphi)-\mathcal{U} _{t} \right \| _{\mathcal{H} }\ge\rho )\\&\ge \mu _{t}(\Gamma _{\rho })-\frac{\widetilde{L _{7}} }{\rho ^{2}} e^{-\eta^{*} t}.
\end{split}
	\end{equation}
 Thus it turns out from \eqref{r37}, \eqref{r38} and Fatou’s Lemma that
\begin{eqnarray*}
\begin{split}
		{\mu }'(\Gamma )&\ge\underset{L \to \infty}{\lim\inf}  \frac{1}{L} \int _{\mathcal{H} }\int_{0}^{L}[\mu _{t}(\Gamma _{\rho })-\frac{\widetilde{L _{7}} }{\rho ^{2}} e^{-\eta^{*} t}]\text{d}t{\mu }'(\text{d}\varphi  )\\&\ge \int _{\mathcal{H} }[\underset{L \to \infty}{\lim\inf} \frac{1}{L} \int_{0}^{L}\mu _{t}(\Gamma _{\rho })\text{d}t]{\mu }'(\text{d}\varphi  )\\&=\mu^{*} (\Gamma _{\rho }).
\end{split}
	\end{eqnarray*}
Let $\rho\to 0$, we have ${\mu }'(\Gamma )\ge \mu ^{*}(\Gamma )$, which implies that ${\mu }'= \mu ^{*}$ by Lemma 2.2.
The proof is complete. \quad $\Box$
~\\

In the following, we will present the uniformly exponential mixing of the  measure $\mu^{*}$ of \eqref{rrr1}.
~\\
\\\textbf{Theorem 5.2.} \emph{Under assumptions \textbf{(A1)}$-$\textbf{(A6)}, as $\lambda<\eta _{1}-\eta _{2}-\frac{145\eta _{3}}{2}  $,  the  measure $\mu^{*}$ of \eqref{rrr1} is uniformly exponential mixing in the sense of Wasserstein metric. More precisely, for any $t\ge0$ and $\nu \in \mathcal{P} (\mathcal{H} )$,
\begin{eqnarray*}
		\left \| P^{*}_{t}\nu -\mu ^{*} \right \|_{BL} \le L_{7}e^{-\frac{\eta^{*} }{2}t }[G^{*}+\int _{\mathcal{H}  }\left \| \phi \right \|^{2}_{\mathcal{H} }\nu (\text{d}\phi )]^{\frac{1}{2} },
	\end{eqnarray*}
where $G^{*}=-\frac{G_{1}+G_{4}}{G_{2}+G_{5}}$.
\\\textbf{proof}} For all $F\in C_{b}(\mathcal{H})$ and $\nu \in \mathcal{P} (\mathcal{H} )$, by \eqref{1}, \eqref{r118} and Chapman–Kolmogorov equation we obtain that
\begin{equation}\label{g10}
\begin{split}
		&\left \| P^{*}_{t}\nu -\mu ^{*} \right \|_{BL} \\&=\left \| P^{*}_{t}\nu -P^{*}_{t}\mu ^{*} \right \|_{BL}\\&=\sup_{\left \| F \right \|_{BL}\le1 } \left | \int_{\mathcal{H} }  F(\psi )\int_{\mathcal{H} }p(t,0,\phi, \text{d}\psi )\nu(\text{d}\phi )- \int_{\mathcal{H} } F(\psi )\int_{\mathcal{H} }p(t,0,\xi, \text{d}\psi )\mu ^{*}(\text{d}\xi ) \right |
\\&=\sup_{\left \| F \right \|_{BL}\le1 } \left |  \int_{\mathcal{H} }\mathbb{E} F(u_{t}(0,\phi ))\nu(\text{d}\phi )- \int _{\mathcal{H} }\mathbb{E} F(u_{t}(0,\xi  )) \mu ^{*}(\text{d}\xi )\right |\\&\le\sup_{\left \| F \right \|_{BL}\le1 }\{\left \| F \right \|_{Lip}(  \int_{\mathcal{H} }\int_{\mathcal{H} }\mathbb{E}\left \| u_{t}(0,\phi )-u_{t}(0,\xi ) \right \| _{\mathcal{H} }^{2}\mu ^{*}(\text{d}\xi)\nu(\text{d}\phi ) )^{\frac{1}{2} }\}\\&\le L_{7}e^{-\frac{\eta^{*} }{2}t }(  \int_{\mathcal{H} }\int_{\mathcal{H} }\left \|  \xi-\phi  \right \|_{\mathcal{H} } ^{2}\mu ^{*}(\text{d}\xi)\nu(\text{d}\phi ) )^{\frac{1}{2} }.
\end{split}
	\end{equation}
In addition, by \eqref{g11}  we have
\begin{equation}\label{g12}
\begin{split}
\int_{\mathcal{H} }\left \|  \xi  \right \|_{\mathcal{H} } ^{2}\mu^{*}(\text{d}\xi)=\lim_{L \to \infty}\frac{1}{L} \int_{0}^{L}  \int_{\mathcal{H} }\left \|  \xi  \right \|_{\mathcal{H} } ^{2}\mu_{t}(\text{d}\xi)\text{d}t=\lim_{L \to \infty}\frac{1}{L} \int_{0}^{L} \mathbb{E} \left \| \mathcal{U} _{t} \right \|^{2}_{\mathcal{H} }  \text{d}t\le-\frac{G_{1}+G_{4}}{G_{2}+G_{5}}.
\end{split}
	\end{equation}
By \eqref{g10} and \eqref{g12} we obtain
\begin{eqnarray*}
		\left \| P^{*}_{t}\nu -\mu^{*}\right \|_{BL} \le L_{7}e^{-\frac{\eta^{*} }{2}t }[-\frac{G_{1}+G_{4}}{G_{2}+G_{5}}+\int _{\mathcal{H}  }\left \| \phi \right \|^{2}_{\mathcal{H} }\nu (\text{d}\phi )]^{\frac{1}{2} }.
	\end{eqnarray*}
By Definition 2.5 of \cite{R7}, we verify that the  measure $\mu^{*}$  is uniformly exponential mixing.
This concludes the proof.\quad $\Box $
~\\
\section{\textup{SLLN and CLT }}
In this section, the strong law of large numbers and the central limit theorem for the solution map of \eqref{rrr1} are established based on uniform exponential mixing Markov processes. Before some details are given, the following preliminaries are introduced.

Let $C_{\mathcal{K} }(\mathcal{H} )$ denote the family of all continuous functionals on $\mathcal{H}$ such that
\begin{eqnarray*}
\begin{split}
		\left \| F \right \| _{\mathcal{K}  }:=\sup_{\varphi \in \mathcal{H} }\frac{|F(\varphi)|}{\mathcal{K}(\left \| \varphi  \right \|_{\mathcal{H} } )  }+\sup_{\varphi_{1}\ne \varphi_{2}}\frac{\left | F(\varphi_{1})-F(\varphi_{2}) \right | }{\left \|\varphi_{1}-\varphi_{2}  \right \|_{\mathcal{H} }[\mathcal{K}(\left \| \varphi_{1}  \right \|_{\mathcal{H} } )+\mathcal{K}(\left \| \varphi_{2}  \right \|_{\mathcal{H} } )]  }<\infty ,
\end{split}
	\end{eqnarray*}
where $\mathcal{K} (t) $ is an increasing continuous functions and $\mathcal{K} (t) >0$ for $t\ge 0$. In this section, we choose the $\mathcal{K} (t) $ as a bounded function, i.e., there exists $M^{\star }>0$ such that $\left | \mathcal{K} (t) \right |\le M^{\star } $ for all $t\ge 0$.

\subsection{Strong law of large numbers}
Firstly, based on the uniformly exponential mixing of the  measure $\mu^{*}$ of \eqref{rrr1}, we prove the strong law of large numbers for  a class of stochastic functional partial differential equations with monotone
coefficients.
~\\
\\\textbf{Theorem 6.1.}\emph{ Under assumptions \textbf{(A1)}$-$\textbf{(A6)}, as $\lambda<\eta _{1}-\eta _{2}-\frac{145\eta _{3}}{2}  $,  for any $\varphi\in \mathcal{H}$ and $F\in C_{\mathcal{K} }(\mathcal{H} ) $, we obtain the following conclusions:
\begin{enumerate}[1)]
        \item There exists a constant $C_{1}$ such that
        \begin{equation}\label{g16}
\begin{split}
		\mathbb{E} \left |\frac{1}{t}\int_{0}^{t}F(u_{s}(\varphi ))\text{d}s -\int _{\mathcal{H} } F(\phi ) \mu ^{*}(\text{d}\phi ) \right | ^{2}\le C_{1}(1+\left \| \varphi  \right \|^{2}_{\mathcal{H}  }  )\left \| F \right \|^{2}_{\mathcal{K}}t^{-1}, \quad t\ge 1;
\end{split}
	\end{equation}
		\item  Let $\alpha >0$ and $\left \lceil \alpha \right \rceil$  denote the integer part of $\alpha$. For any fix $\left \lceil \alpha \right \rceil$, we obtain that for any $\varepsilon \in (0,\frac{1}{2(\left \lceil \alpha \right \rceil+2)}  )$, there exists a constant $C_{2}>0$ such that
\begin{equation}\label{g17}
\begin{split}
		\left |\frac{1}{t}\int_{0}^{t}F(u_{s}(\varphi ))\text{d}s -\int _{\mathcal{H} } F(\phi ) \mu ^{*}(\text{d}\phi ) \right |\le C_{2}\left \| F \right \|_{\mathcal{K}}t^{-\frac{1}{2(\left \lceil \alpha \right \rceil+2)} +\varepsilon }, \quad t\ge T_{\varepsilon}(\omega ),\quad \mathbb{P} -a.s.,
\end{split}
	\end{equation}
where  the random time $T_{\varepsilon}(\omega )$ is $\mathbb{P}-$a.s. finite. Moreover, for any $ {\alpha }' \in (0,\frac{\alpha }{\left \lceil \alpha \right \rceil+2} )$, there exists a constant ${C}'_{1}$ such that
\begin{equation}\label{g18}
\begin{split}
		\mathbb{E} T_{\varepsilon}^{{\alpha }' }(\omega )\le 1+ \frac{{C}'_{1}}{\alpha -(\left \lceil \alpha \right \rceil+2){\alpha }' } (1+\left \| \varphi  \right \|^{2}_{\mathcal{H}  }  )\left \| F \right \|^{2}_{\mathcal{K}}.
\end{split}
	\end{equation}
	\end{enumerate}
\textbf{proof}}
\begin{enumerate}[1)]
\item For any given $F\in C_{\mathcal{K} }(\mathcal{H} )$ and $\varphi \in \mathcal{H} $, it follows from Theorem 5.2 and Chapman-Kolmogorov equation that
\begin{equation}\label{g15}
\begin{split}
\left | P_{t}F(\varphi )- \int _{\mathcal{H} }F(\phi )\mu ^{*}(\text{d}\phi )\right | &=\left |\int _{\mathcal{H} } F(\phi  )p(t,0,\varphi ,\text{d}\phi )- \int _{\mathcal{H} }F(\phi )\mu ^{*}(\text{d}\phi )\right |\\&=\left | \int _{\mathcal{H} }\int _{\mathcal{H} }F(\phi  )p(t,0,\xi  ,\text{d}\phi )p(0,0,\varphi ,\text{d}\xi  )- \int _{\mathcal{H} }F(\phi )\mu ^{*}(\text{d}\phi )\right |\\&=\left | \int _{\mathcal{H} }F(\phi  )P^{*}_{t}p(0,0,\varphi ,\text{d}\phi )- \int _{\mathcal{H} }F(\phi )\mu ^{*}(\text{d}\phi )\right |\\&\le \left \| F \right \| _{BL}\cdot \left \| P^{*}_{t}\mu _{0}^{\varphi } -\mu ^{*}\right \|_{BL}\\&\le 2L_{7}M^{\star }\left \| F \right \| _{\mathcal{K}}e^{-\frac{\eta^{*} }{2}t }[G^{*}+\int _{\mathcal{H}  }\left \| \phi \right \|^{2}_{\mathcal{H} }p(0,0,\varphi ,\text{d}\phi )]^{\frac{1}{2} }\\&=2L_{7}M^{\star }\left \| F \right \| _{\mathcal{K}}e^{-\frac{\eta^{*} }{2}t }[G^{*}+\left \| \varphi\right \|^{2}_{\mathcal{H} }]^{\frac{1}{2} },
\end{split}
	\end{equation}
where $\mu _{0}^{\varphi }(\Gamma)=p(0,0,\varphi,\Gamma)$ for any $\Gamma\in \mathcal{B}(\mathcal{H} )$. Hence  the assumptions in Lemma 2.1 of \cite{R8} holds for $\mathbb{B}=\mathcal{H}$, $k=1$, $\varphi (t)=2L_{7}M^{\star }e^{-\frac{\eta^{*} }{2}t }$ and $\psi (\left \| \varphi  \right \|_{\mathcal{H} } )=[G^{*}+\left \| \varphi \right \|_{\mathcal{H} }^{2}]^{\frac{1}{2} }$, then the desired conclusion \eqref{g16} holds(see also  Proposition 2.6 of \cite{R7}).
\item Step 1:  There is no loss of generality in assuming that $\int _{\mathcal{H} } F(\phi ) \mu ^{*}(\text{d}\phi )=0$ and $\left \| F \right \| _{\mathcal{K} }\le 1$.

For $k\ge1$, let
\begin{eqnarray*}
t_{k}=k^{\left \lceil \alpha \right \rceil+2}, \quad \upsilon =\frac{\left \lceil \alpha \right \rceil+1-\alpha}{2},
	\end{eqnarray*}
 and  consider the events
\begin{eqnarray*}
\mathfrak{D}(k) =\left \{ \omega \in \Omega :\left | \frac{1}{t_{k}} \int_{0}^{t_{k}}F(u_{s}(\varphi ))\text{d}s  \right |>\frac{1}{k^{\upsilon }}   \right \}.
	\end{eqnarray*}
By Chebyshev inequality and \eqref{g16}, for any $\varphi \in \mathcal{H} $ and $k\ge1$, we obtain
\begin{equation}\label{g19}
\begin{split}
		\mathbb{P} (\mathfrak{D}(k) )\le \mathbb{E} \left | \frac{1}{t_{k}} \int_{0}^{t_{k}}F(u_{s}(\varphi ))\text{d}s  \right |^{2} k^{2\upsilon }\le C_{1}(1+\left \| \varphi  \right \|^{2}_{\mathcal{H}  }  )\left \| F \right \|^{2}_{\mathcal{K}}k^{-1-\alpha }.
\end{split}
	\end{equation}
We shall assume that $k^{*}(\omega )\ge0$ is the smallest integer such that
\begin{equation}\label{g35}
\left | \frac{1}{t_{k}} \int_{0}^{t_{k}}F(u_{s}(\varphi ))\text{d}s  \right |\le k^{-\upsilon}=t_{k}^{- \frac{\upsilon}{\left \lceil \alpha \right \rceil+2} },
	\end{equation}
which  implies  that the random variable $k^{*}(\omega )$ is $\mathbb{P}-$a.s. finite by Borel–Cantelli lemma and  \eqref{g19}. For any $k\ge k^{*}(\omega )$ and $t\in (t_{k},t_{k+1})$, we obtain
\begin{equation}\label{g38}
\begin{split}
		\left |\frac{1}{t}\int_{0}^{t}F(u_{s}(\varphi ))\text{d}s  \right |&\le \left | \frac{1}{t} \int_{0}^{t}F(u_{s}(\varphi ))\text{d}s- \frac{1}{t_{k+1}} \int_{0}^{t_{k+1}}F(u_{s}(\varphi ))\text{d}s \right |\\&+\left | \frac{1}{t_{k+1}} \int_{0}^{t_{k+1}}F(u_{s}(\varphi ))\text{d}s \right |.
\end{split}
	\end{equation}
Clearly, it's not hard to get
\begin{equation}\label{g32}
\begin{split}
		&\left | \frac{1}{t} \int_{0}^{t}F(u_{s}(\varphi ))\text{d}s- \frac{1}{t_{k+1}} \int_{0}^{t_{k+1}}F(u_{s}(\varphi ))\text{d}s \right |\\&\le\frac{1}{t}\left | \int_{0}^{t}F(u_{s}(\varphi ))\text{d}s- \int_{0}^{t_{k+1}}F(u_{s}(\varphi ))\text{d}s \right |+(\frac{1}{t_{k}} -\frac{1}{t_{k+1}} ) \left |\int_{0}^{t_{k+1}}F(u_{s}(\varphi ))\text{d}s  \right |.
\end{split}
	\end{equation}
By the definition of norm $\left \| \cdot  \right \|  _{\mathcal{K} }$ and the boundedness of function $\mathcal{K}$, we derive
\begin{equation}\label{g33}
\begin{split}
	\frac{1}{t}\left | \int_{0}^{t}F(u_{s}(\varphi ))\text{d}s- \int_{0}^{t_{k+1}}F(u_{s}(\varphi ))\text{d}s \right |\le \frac{1}{t} \int_{t}^{t_{k+1}} \mathcal{K} (\left \| u_{s}(\varphi ) \right \| )\text{d}s\le M^{*}\frac{t_{k+1}-t_{k}}{t_{k}},
\end{split}
	\end{equation}
and by \eqref{g35}, we have
\begin{equation}\label{g36}
\begin{split}
	(\frac{1}{t_{k}} -\frac{1}{t_{k+1}} ) \left |\int_{0}^{t_{k+1}}F(u_{s}(\varphi ))\text{d}s\right |\le\frac{t_{k+1}-t_{k}}{t_{k}}t_{k+1}^{- \frac{\upsilon}{\left \lceil \alpha \right \rceil+2} }.
\end{split}
	\end{equation}
In addition, there exists a constant $C=C(\left \lceil \alpha \right \rceil)>0$ such that
\begin{equation}\label{g37}
\begin{split}
	\frac{t_{k+1}-t_{k}}{t_{k}}\le C\frac{1}{k+1} =Ct_{k+1}^{-\frac{1}{\left \lceil \alpha \right \rceil+2} }.
\end{split}
	\end{equation}
 Substituting \eqref{g32}-\eqref{g37} into \eqref{g38} gives that for any $k\ge k^{*}(\omega )$ and $t\in (t_{k},t_{k+1})$,
\begin{eqnarray*}
\begin{split}
		\left |\frac{1}{t}\int_{0}^{t}F(u_{s}(\varphi ))\text{d}s  \right |&\le CM^{*}t_{k+1}^{-\frac{1}{\left \lceil \alpha \right \rceil+2} }+Ct_{k+1}^{-\frac{1}{\left \lceil \alpha \right \rceil+2} }+t_{k+1}^{-\frac{1}{2(\left \lceil \alpha \right \rceil+2)} +\frac{\alpha -\left \lceil \alpha \right \rceil}{2(\left \lceil \alpha \right \rceil+2)}}\\
&\le C_{2}t^{-\frac{1}{2(\left \lceil \alpha \right \rceil+2)} +\frac{\alpha -\left \lceil \alpha \right \rceil}{2(\left \lceil \alpha \right \rceil+2)}}.
\end{split}
	\end{eqnarray*}
We fix arbitrary $\varepsilon \in (0,\frac{1}{2(\left \lceil \alpha \right \rceil+2)})$ and $F\in C_{\mathcal{K} }(\mathcal{H} ) $. Let $\alpha =2\varepsilon (\left \lceil \alpha  \right \rceil +2)+\left \lceil \alpha  \right \rceil$, which implies that \eqref{g17} holds.

Step 2: In the proof in Step 1, we get that \eqref{g17} holds with $T_{\varepsilon}(\omega )=t_{k^{*}(\omega )}=[k^{*}(\omega )]^{\left \lceil \alpha \right \rceil+2}$. In the following, we prove \eqref{g18}.
For any ${\alpha }' \in (0,\frac{\alpha }{\left \lceil \alpha \right \rceil+2} )$ we have
\begin{eqnarray*}
\begin{split}
		\mathbb{E} [T_{\varepsilon}(\omega )]^{{\alpha }' }&=\mathbb{E}[k^{*}(\omega )]^{(\left \lceil \alpha \right \rceil+2){\alpha }' }\\&\le\sum_{k=1}^{\infty }\mathbb{P}(k^{*}=k) k^{(\left \lceil \alpha \right \rceil+2){\alpha }' }\\&\le1+ \sum_{k=2}^{\infty }\mathbb{P}(\mathfrak{D}(k) ) k^{(\left \lceil \alpha \right \rceil+2){\alpha }' }\\&\le1+C_{1}(1+\left \| \varphi  \right \|^{2}_{\mathcal{H}  }  )\left \| F \right \|^{2}_{\mathcal{K}}\sum_{k=2}^{\infty }k^{-1-\alpha+(\left \lceil \alpha \right \rceil+2){\alpha }' }\\&\le 1+\frac{{C}'_{1} }{\alpha -(\left \lceil \alpha \right \rceil+2){\alpha }' }  (1+\left \| \varphi  \right \|^{2}_{\mathcal{H}  }  )\left \| F \right \|^{2}_{\mathcal{K}}.
\end{split}
	\end{eqnarray*}
	\end{enumerate}
The proof of Theorem 6.1 is complete. $\quad \Box$
\subsection{Central limit theorem}
To state the CLT, we first introduce the corrector $\Upsilon :C_{\mathcal{K} }(\mathcal{H} )\longrightarrow \mathbb{R} $ such that
\begin{eqnarray*}
\Upsilon [F(\varphi )] =\int_{0}^{\infty } [P_{t}F(\varphi )-\int _{\mathcal{H}}F(\phi )\mu ^{*}(\text{d}\phi ) ]\text{dt},
	\end{eqnarray*}
for any $F\in C_{\mathcal{K} }(\mathcal{H} )$ and $\varphi\in\mathcal{H} $. In this subsection, let us fix $\varphi\in\mathcal{H} $ and an arbitrary function $F\in C_{\mathcal{K} }(\mathcal{H} )$ such that $\int _{\mathcal{H}}F(\phi )\mu ^{*}(\text{d}\phi )=0$ and set
\begin{eqnarray*}
\mathcal{M} _{t}^{F}=\int_{0}^{\infty } [\mathbb{E} (F(u_{r}(\varphi ))|\mathscr{F}_{t} )-\mathbb{E}(F(u_{r}(\varphi ))|\mathscr{F}_{0} ) ]\text{d}r,
\end{eqnarray*}
\begin{eqnarray*}
\mathcal{S} [F(\varphi )]=\mathbb{E} \left | \int_{0}^{1}F(u_{r}(\varphi ))\text{d}t+\Upsilon [F(u_{1}(\varphi ))]-\Upsilon [F(\varphi )]  \right |^{2}.
\end{eqnarray*}
By \eqref{g15} we have
\begin{equation}\label{g20}
\begin{split}
		\left | P_{t}F(\varphi ) \right |  &\le 2L_{7}M^{\star }\left \| F \right \| _{\mathcal{K}}e^{-\frac{\eta^{*} }{2}t }[G^{*}+\left \| \varphi\right \|^{2}_{\mathcal{H} }]^{\frac{1}{2} }\\&=M^{\star }_{\ast }\left \| F \right \| _{\mathcal{K}}e^{-\frac{\eta^{*} }{2}t }(1+\left \| \varphi\right \|^{2}_{\mathcal{H} })^{\frac{1}{2} },
\end{split}
	\end{equation}
which implies that $\Upsilon [F(\varphi )] =\int_{0}^{\infty } P_{t}F(\varphi )\text{dt}$ and $\mathcal{S} [F(\varphi )]$ are well-defined by \eqref{g20} and \eqref{26}.
~\\
\\\textbf{Lemma 6.2.} \emph{Assume the conditions of Theorem 5.2 hold. Then  $\mathcal{M} _{t}^{F}$ is a well-defined square integrable martingale.
\\\textbf{proof}} Firstly, it is proved that $\mathcal{M} _{t}^{F}$ is a martingale. Indeed, by  the dominated convergence theorem, \eqref{g20} and \eqref{26}, for any $t>s\ge0$, we obtain
\begin{eqnarray*}
\begin{split}
\mathbb{E} (\mathcal{M}_{t}^{F}|\mathscr{F}_{s}  )&=\int_{0}^{\infty } [\mathbb{E}( \mathbb{E} (F(u_{r}(\varphi ))|\mathscr{F}_{t} )|\mathscr{F}_{s}) -\mathbb{E} (\mathbb{E}(F(u_{r}(\varphi ))|\mathscr{F}_{0} )|\mathscr{F}_{s})  ]\text{d}r\\
&=\int_{0}^{\infty } [\mathbb{E} (F(u_{r}(\varphi ))|\mathscr{F}_{s} )-\mathbb{E}(F(u_{r}(\varphi ))|\mathscr{F}_{0} ) ]\text{d}r\\
&=\mathcal{M}_{s}^{F}.
\end{split}
\end{eqnarray*}
In addition, it's not hard to get that $\mathbb{E} \mathcal{M}_{t}^{F}=0$.

Next, by \eqref{r118} and the  Markov properties of solution map $u_{t}$, we obtain
\begin{equation}\label{g21}
\begin{split}
\mathcal{M} _{t}^{F}&=\int_{0}^{\infty } [\mathbb{E} (F(u_{r}(\varphi ))|\mathscr{F}_{t} )-\mathbb{E}(F(u_{r}(\varphi ))|\mathscr{F}_{0} ) ]\text{d}r\\
&=\int_{0}^{t }F(u_{r}(\varphi )) \text{d}r+\int_{t}^{\infty } \mathbb{E} (F(u_{r}(\varphi ))|\mathscr{F}_{t} )\text{d}r-\int_{0}^{\infty }\mathbb{E}(F(u_{r}(\varphi )) ) \text{d}r\\
&=\int_{0}^{t }F(u_{r}(\varphi ))-P_{r}F(\varphi) \text{d}r+\int_{t}^{\infty } P_{r-t}F(u_{t}(\varphi ))-P_{r}F(\varphi)\text{d}r.
\end{split}
\end{equation}
Thus, \eqref{g20}, \eqref{26} and the definition of $\left \| \cdot  \right \| _{\mathcal{K} }$ yield
\begin{eqnarray*}
\begin{aligned}
&\mathbb{E} \left | \mathcal{M} _{t}^{F} \right | ^{2}\\&
=\mathbb{E} \left |\int_{0}^{t }F(u_{r}(\varphi ))-P_{r}F(\varphi) \text{d}r+\int_{t}^{\infty } P_{r-t}F(u_{t}(\varphi ))-P_{r}F(\varphi)\text{d}r\right | ^{2}\\
&=\mathbb{E} \left | \left \| F \right \|_{\mathcal{K} }\{\int_{0}^{t}\mathcal{K}(u_{r}(\varphi ))\text{d}r+M^{\star }_{\ast }[\int_{t}^{\infty }e^{-\frac{\eta^{*} }{2}(r-t) }(1+\left \| u_{t}(\varphi ) \right \|_{\mathcal{H} }^{2} )^{\frac{1}{2} } \text{d}r+\int_{0}^{\infty }e^{-\frac{\eta^{*} }{2}r }(1+\left \| \varphi  \right \|_{\mathcal{H} }^{2} )^{\frac{1}{2} } \text{d}r] \}  \right |^{2}\\
&\le M^{\star }_{\ast }(t) \left \| F \right \|_{\mathcal{K} }[1+\mathbb{E}\left \| u_{t}(\varphi ) \right \|_{\mathcal{H} }^{2}+\left \| \varphi  \right \|_{\mathcal{H} }^{2} ],
\end{aligned}
\end{eqnarray*}
where $M^{\star }_{\ast }(t)$ is an increasing function, which implies that $\mathbb{E} \left | \mathcal{M} _{t}^{F} \right | ^{2}< \infty  $ for $t\in \mathbb{R} $.\quad $\Box$
~\\

For any integer $k\ge 1$, by \eqref{g21}, we have
\begin{equation}\label{g22}
\begin{aligned}
\mathcal{M} _{k}^{F}=\mathcal{M} _{k-1}^{F}+\int_{k-1}^{k} F(u_{r}(\varphi ))\text{d}r+\int_{0}^{\infty } P_{r}F(u_{k}(\varphi ))-P_{r}F(u_{k-1}(\varphi ))\text{d}r.
\end{aligned}
\end{equation}
Next consider the conditional variance for $\{\mathcal{M}_{n}^{F} \}_{n\in\mathbb{Z}_{+}    }$,
\begin{eqnarray*}
\begin{aligned}
\Xi ^{2}_{n}( \mathcal{M}^{F}) =\sum_{i=1}^{n} \mathbb{E} ((\mathcal{M}_{i}^{F} -\mathcal{M}_{i-1}^{F})^{2}|\mathscr{F}_{i-1} ).
\end{aligned}
\end{eqnarray*}
By \eqref{g22} and the Markov properties of  $u_{t}$, we obtain
\begin{eqnarray*}
\begin{aligned}
\mathbb{E} ((\mathcal{M}_{i}^{F} -\mathcal{M}_{i-1}^{F})^{2}|\mathscr{F}_{i-1} )&=\mathbb{E} (\left | \int_{i-1}^{i} F(u_{r}(\varphi ))\text{d}r+\int_{0}^{\infty } P_{r}F(u_{i}(\varphi ))-P_{r}F(u_{i-1}(\varphi ))\text{d}r \right |^{2} |\mathscr{F}_{i-1}) \\&
=\mathcal{S} [F(u_{i-1}(\varphi) )],
\end{aligned}
\end{eqnarray*}
which implies
\begin{equation}\label{g23}
\begin{aligned}
\Xi ^{2}_{n}( \mathcal{M}^{F}) =\sum_{i=0}^{n-1}\mathcal{S} [F(u_{i}(\varphi) )].
\end{aligned}
\end{equation}
~\\
\\\textbf{Lemma 6.3.} \emph{Assume the conditions of Theorem 5.2 hold. Then there exist constants $c_{1},c_{2}$ such that
\begin{eqnarray*}
\begin{aligned}
 \mathbb{E} (\sup_{t\in (n,n+1)}e^{c_{1}\left \| u_{t}(\varphi ) \right \|^{2}_{\mathcal{H} } } )\le e^{c_{2}(1+\left \| \varphi  \right \|^{2}_{\mathcal{H} } )},
\end{aligned}
\end{eqnarray*}
for any $n\ge0$ and $\varphi \in \mathcal{H} $.
\\\textbf{proof} }Similar to the proof of Lemma 2.1 in \cite{R11} and Lemma 4.3 in \cite{R8}.
~\\
\\\textbf{Lemma 6.4.}\emph{ Assume the conditions of Theorem 5.2 hold. Then for any $F\in C_{\mathcal{K} }(\mathcal{H} )$ with $\int _{\mathcal{H}}F(\phi )\mu ^{*}(\text{d}\phi )=0$, we obtain
\begin{eqnarray*}
\begin{aligned}
 0\le \int _{\mathcal{H}}\mathcal{S}  [F(\phi )]\mu ^{*}(\text{d}\phi ) =2 \int _{\mathcal{H}}F(\phi )\Upsilon[F(\phi )]\mu ^{*}(\text{d}\phi )<\infty.
\end{aligned}
\end{eqnarray*}
\textbf{proof}} Similar to the proof of  Lemma 4.1 in \cite{R8}.
~\\
\\\textbf{Remark 6.5.} It follows from \eqref{g10} that for any $\nu_{1},\nu_{2}\in\mathcal{P}  (\mathcal{H} )$
\begin{equation}\label{g25}
\begin{aligned}
\left \| P^{*}_{t}\nu_{1} -P^{*}_{t}\nu_{2} \right \|_{BL} \le L_{7}e^{-\frac{\eta^{*} }{2}t }(  \int_{\mathcal{H} }\int_{\mathcal{H} }\left \|  \xi-\phi  \right \|_{\mathcal{H} } ^{2}\nu_{1}(\text{d}\xi )\nu_{2}(\text{d}\phi ))^{\frac{1}{2} }.
\end{aligned}
\end{equation}
In addition, by \eqref{g15} and \eqref{g25},  we have
\begin{equation}\label{g26}
\begin{aligned}
\left | P_{t}F(\varphi )- P_{t}F(\xi )\right |&\le \left \| F \right \| _{BL}\cdot \left \| P^{*}_{t}\mu _{0}^{\varphi } -P^{*}_{t}\mu _{0}^{\xi }\right \|_{BL}\\&\le 2L_{7}M^{\star }\left \| F \right \| _{\mathcal{K}}e^{-\frac{\eta^{*} }{2}t }[\int _{\mathcal{H}  }\int _{\mathcal{H}  }\left \| \phi_{1}-\phi_{2} \right \|^{2}_{\mathcal{H} }p(0,0,\varphi ,\text{d}\phi_{1} )p(0,0, \xi,\text{d}\phi_{2} )]^{\frac{1}{2} }\\
&=2L_{7}M^{\star }\left \| F \right \| _{\mathcal{K}}e^{-\frac{\eta^{*} }{2}t }[\int _{\mathcal{H}  }\left \| \varphi-\phi_{2} \right \|^{2}_{\mathcal{H} }p(0,0,\xi ,\text{d}\phi_{2} )]^{\frac{1}{2} }\\
&=2L_{7}M^{\star }\left \| F \right \| _{\mathcal{K}}e^{-\frac{\eta^{*} }{2}t }\left \| \varphi-\xi \right \|_{\mathcal{H}},
\end{aligned}
\end{equation}
where $\mu _{0}^{\varphi }(\Gamma)=p(0,0,\varphi,\Gamma)$ and $\mu _{0}^{\xi }(\Gamma)=p(0,0,\xi,\Gamma)$. Thus, by \eqref{g25} and \eqref{g26}, it is similar to  Lemma 4.2 in \cite{R8},  and we can obtain that there exists a constant $\overline{M} >0$ such that for any $F\in C_{\mathcal{K} }(\mathcal{H} )$ and $\varphi \in \mathcal{H} $
\begin{equation}\label{g27}
\begin{aligned}
 \left \| \mathcal{S}[F]  \right \| _{\mathcal{K} }\le \overline{M} \left \| F \right \| _{\mathcal{K} }^{2}.
\end{aligned}
\end{equation}

For any the non-negative constant $\zeta $, let $\Pi _{\zeta }(\cdot )$ be the normal distribution function with zero mean and variance $\zeta$, where $$\Pi _{0 }(y)=\left\{\begin{matrix}
  1,& y\ge0,\\
  0,& y<0.
\end{matrix}\right.$$

Then based on the above preparation, we will give the central limit theorem below.
~\\
\\\textbf{Theorem 6.6.}\emph{ Under assumptions \textbf{(A1)}$-$\textbf{(A6)} and $\lambda<\eta _{1}-\eta _{2}-\frac{145\eta _{3}}{2}  $,  for any $\varphi\in \mathcal{H}$ and $F\in C_{\mathcal{K} }(\mathcal{H} ) $ with $\int _{\mathcal{H}}F(\phi )\mu ^{*}(\text{d}\phi )=0$, let $\zeta =(\int _{\mathcal{H}}\mathcal{S}  [F(\phi )]\mu ^{*}(\text{d}\phi ))^{\frac{1}{2} }\in [0,\infty )$ and we have the following conclusions:
\begin{enumerate}[(\textbf{a})]
		\item When $\zeta>0$, for $\epsilon  \in [0,\frac{1}{5} )$, there exists an increasing function $\mathcal{G} _{\epsilon  }:\mathbb{R}_{+}\times \mathbb{R} _{+}\to \mathbb{R} _{+}$ such that
\begin{eqnarray*}
\begin{aligned}
 \sup_{y\in \mathbb{R} } \left | \mathbb{P} (\frac{1}{\sqrt{t} }\int_{0}^{t}F(u_{r}(\varphi ))\text{d}r\le y  )-\Pi _{\zeta } (y)\right | \le \mathcal{G}_{\epsilon  }(\left \| F \right \| _{\mathcal{K} },\left \| \varphi  \right \|_{\mathcal{H} } )t^{-\frac{1}{5}+\epsilon   },
\end{aligned}
\end{eqnarray*}
for any $\varphi \in \mathcal{H} $ and $t>0$;
\end{enumerate}
\begin{enumerate}[(\textbf{b})]
        \item When $\zeta=0$, there exists an increasing function $\mathcal{G} :\mathbb{R} _{+}\times \mathbb{R} _{+}\to \mathbb{R} _{+}$ such that
\begin{eqnarray*}
\begin{aligned}
 \sup_{y\in \mathbb{R} }[ (\left | y \right |\wedge 1 ) \left | \mathbb{P} (\frac{1}{\sqrt{t} }\int_{0}^{t}F(u_{r}(\varphi ))\text{d}r\le y  )-\Pi _{0 } (y)\right | ]\le \mathcal{G}(\left \| F \right \| _{\mathcal{K} },\left \| \varphi  \right \|_{\mathcal{H} } )t^{-\frac{1}{4} },
\end{aligned}
\end{eqnarray*}
for any $\varphi \in \mathcal{H} $ and $t>0$.
	\end{enumerate}
\textbf{Proof of (a)}} Theorem 5.2 and Lemma 6.3 imply that the assumptions of Theorem 2.8 in \cite{R7} hold. Thus,  we can obtain that for $\varepsilon \in [\frac{1}{20},\frac{1}{4} )$ and $q=1$, there exists an increasing continuous function $h _{\varepsilon }:\mathbb{R} _{+}\times \mathbb{R} _{+}\to \mathbb{R} _{+}$ such that
\begin{equation}\label{g28}
\begin{aligned}
 &\sup_{y\in \mathbb{R} } \left | \mathbb{P} (\frac{1}{\sqrt{t} }\int_{0}^{t}F(u_{r}(\varphi ))\text{d}r\le y  )-\Pi _{\zeta } (y)\right |\\ &\le t^{-\frac{1}{4}+\varepsilon  }h_{\varepsilon }(\left \| F \right \| _{\mathcal{K} },\left \| \varphi  \right \|_{\mathcal{H} } )+\zeta^{-4}\left \lceil t \right \rceil ^{(1-4\varepsilon)}\mathbb{E} \left | \frac{\Xi ^{2}_{\left \lceil t \right \rceil}( \mathcal{M}^{F})}{\left \lceil t \right \rceil } -\zeta ^{2} \right |^{2}.
\end{aligned}
\end{equation}
 Similar to \eqref{g15}, we have
\begin{eqnarray*}
\begin{aligned}
 \left | P_{t}\mathcal{S}[F(\varphi )] -\zeta^{2} \right | &=\left | P_{t}\mathcal{S}[F(\varphi )] -\int _{\mathcal{H}}\mathcal{S}  [F(\phi )]\mu ^{*}(\text{d}\phi ) \right |\\
 &\le M^{\star }_{\ast }\left \| \mathcal{S}[F] \right \| _{\mathcal{K}}e^{-\frac{\eta^{*} }{2}t }(1+\left \| \varphi\right \|^{2}_{\mathcal{H} })^{\frac{1}{2} },
\end{aligned}
\end{eqnarray*}
which implies that there exists a constant $C_{3}$ such that by   Lemma 2.1 of \cite{R8}
\begin{eqnarray*}
\begin{aligned}
 \mathbb{E} \left | \frac{1}{n} \sum_{i=0}^{n-1}\mathcal{S} [F(u_{i}(\varphi) )]-\zeta^{2} \right |^{2} \le C_{3}n^{-1}\left \| \mathcal{S}[F] \right \|^{2} _{\mathcal{K}}(1+\left \| \varphi\right \|^{2}_{\mathcal{H} }).
\end{aligned}
\end{eqnarray*}
Hence by \eqref{g23}, we arrive at
\begin{equation}\label{g29}
\begin{aligned}
 \mathbb{E} \left |\frac{\Xi ^{2}_{\left \lceil t \right \rceil}( \mathcal{M}^{F})}{\left \lceil t \right \rceil} -\zeta^{2} \right |^{2} \le C_{3}\left \lceil t \right \rceil^{-1}\left \| \mathcal{S}[F] \right \|^{2} _{\mathcal{K}}(1+\left \| \varphi\right \|^{2}_{\mathcal{H} }).
\end{aligned}
\end{equation}
When $\varepsilon\ge \frac{1}{20}$, we have $\left \lceil t \right \rceil^{-4\varepsilon } \le 2^{4\varepsilon }\cdot t^{-\frac{1}{4}+\varepsilon  }$, which implies $\mathcal{G} _{\epsilon  }(\left \| F \right \| _{\mathcal{K} },\left \| \varphi  \right \|_{\mathcal{H} } )=h_{\varepsilon }(\left \| F \right \| _{\mathcal{K} },\left \| \varphi  \right \|_{\mathcal{H} } )+2^{4\varepsilon}C_{3}\overline{M} ^{2}\zeta^{-4}(1+\left \| \varphi\right \|^{2}_{\mathcal{H} }) \left \| F \right \| _{\mathcal{K} }^{4}$ by \eqref{g28}, \eqref{g29} and Remark 6.5. Finally, let $\epsilon =\varepsilon -\frac{1}{20} $ and we  prove (a) of Theorem 6.6.
~\\
\\\textbf{\emph{Proof of (b)}} When $\zeta=0$, we can obtain  from Theorem 2.8 in \cite{R7} that there exists an increasing continuous function $h :\mathbb{R} _{+}\times \mathbb{R} _{+}\to \mathbb{R} _{+}$ such that
\begin{equation}\label{g30}
\begin{aligned}
 &\sup_{y\in \mathbb{R} } \left | \mathbb{P} (\frac{1}{\sqrt{t} }\int_{0}^{t}F(u_{r}(\varphi ))\text{d}r\le y  )-\Pi _{0 } (y)\right |\\ &\le t^{-\frac{1}{4} } h(\left \| F \right \| _{\mathcal{K} },\left \| \varphi  \right \|_{\mathcal{H} } )+\left \lceil t \right \rceil ^{-\frac{1}{2} } \left |\mathbb{E} \Xi ^{2}_{\left \lceil t \right \rceil}( \mathcal{M}^{F}) \right |^{\frac{1}{2}}\\
 &\le t^{-\frac{1}{4} } h(\left \| F \right \| _{\mathcal{K} },\left \| \varphi  \right \|_{\mathcal{H} } )+\left \lceil t \right \rceil ^{-\frac{1}{4} }\left \| \mathcal{S}[F] \right \| _{\mathcal{K}}[C_{3}(1+\left \| \varphi\right \|^{2}_{\mathcal{H} }) ]^{\frac{1}{2} }\\
 &\le\mathcal{G}(\left \| F \right \| _{\mathcal{K} },\left \| \varphi  \right \|_{\mathcal{H} } )t^{-\frac{1}{4} },
\end{aligned}
\end{equation}
where $\mathcal{G}(\left \| F \right \| _{\mathcal{K} },\left \| \varphi  \right \|_{\mathcal{H} } )=h(\left \| F \right \| _{\mathcal{K} },\left \| \varphi  \right \|_{\mathcal{H} } )+16\overline{M} [C_{3}(1+\left \| \varphi\right \|^{2}_{\mathcal{H} }) ]^{\frac{1}{2} }\left \| F \right \| _{\mathcal{K} }^{2}$. The proof of Theorem 6.6 is complete. \quad $\Box$
\section{\textup{Applications}}
In this section, in order to illustrate the validity of our main results,  we give the application to  stochastic generalized porous media equations. It is worth noting that in the example, we mainly consider additive or linear multiplicative noise. Let $\Lambda \subset \mathbb{R}^{n}(n \in \mathbb{N}) $ be an open bounded subset and let $ - \Delta$ have the Dirichlet boundary condition.
~\\
\\\textbf{Example} Consider the stochastic generalized porous media equation:
\begin{equation}\label{g31}
\begin{cases}
\begin{aligned}
 \text{d}u&=[\Delta(\left | u \right | ^{q-2}+a_{1}u) +f(t,u_{t})]\text{d}t+g(t,u_{t})\text{d}W(t),
 \end{aligned}\\
 u_{0}=\varphi\in \mathcal{H},
\end{cases}.
\end{equation}
where $a_{1}\ge0$ and $W$ is a one-dimensional two-sided  cylindrical $Q$-Wiener process with $Q = I$ on $L^{q}(\Lambda )(q>2)$.

Assume that $f$ and $g$ satisfy (\textbf{A5}) with the same constants and have some recurrence. Then we have the following theorem.
~\\
\\\textbf{Theorem 7.1.} \emph{Assume that $a_{1}+\eta _{1}-\eta _{2}-\frac{145\eta _{3}}{2}>0  $, then  we have the
following conclusions:
\begin{enumerate}[(1)]
        \item There exists a unique entrance measure $\mu _{t}\in \mathcal{R}_{2}$ with the same recurrence properties as the coefficients $f$ and $g$ in the sense of Wasserstein metric for  system \eqref{g31}.
		\item  When \eqref{g31} is an autonomous system, i.e., $f(t,\varphi)=f(\varphi)$ and $g(t,\varphi)=g(\varphi)$:
\begin{enumerate}[(i)]
\item There exists a unique  measure $\mu^{*}$ to  system \eqref{g31} which  is uniformly exponentially mixing in the sense of Wasserstein metric, i.e., for any $t\ge0$ and $\nu \in \mathcal{P}$, there exist constants $L, \eta >0$ such that
\begin{eqnarray*}
		\left \| P^{*}_{t}\nu -\mu ^{*} \right \|_{BL} \le Le^{-\eta t }[1+\int _{\mathcal{H}  }\left \| \phi \right \|^{2}_{\mathcal{H} }\nu (\text{d}\phi )]^{\frac{1}{2} };
	\end{eqnarray*}
\item For any  $\varphi\in \mathcal{H}$  and $F\in C_{\mathcal{K} }(\mathcal{H} ) $,    system \eqref{g31} has a unique global solution map $u_{t}(0,\varphi )=u_{t}(\varphi )$ such that
\begin{enumerate}[(\textbf{a})]
        \item there exists a constant $C$ such that
        \begin{eqnarray*}
\begin{split}
		\mathbb{E} \left |\frac{1}{t}\int_{0}^{t}F(u_{s}(\varphi ))\text{d}s -\int _{\mathcal{H} } F(\phi ) \mu ^{*}(\text{d}\phi ) \right | ^{2}\le C(1+\left \| \varphi  \right \|^{2}_{\mathcal{H}  }  )\left \| F \right \|^{2}_{\mathcal{K}}t^{-1}, \quad t\ge 1;
\end{split}
	\end{eqnarray*}
\end{enumerate}
\begin{enumerate}[(\textbf{b})]
		\item   For any fix $\left \lceil \alpha \right \rceil$, we obtain that for any $\varepsilon \in (0,\frac{1}{2(\left \lceil \alpha \right \rceil+2)}  )$, there exists a constant $C>0$ such that
\begin{eqnarray*}
\begin{split}
		\left |\frac{1}{t}\int_{0}^{t}F(u_{s}(\varphi ))\text{d}s -\int _{\mathcal{H} } F(\phi ) \mu ^{*}(\text{d}\phi ) \right |\le C\left \| F \right \|_{\mathcal{K}}t^{-\frac{1}{2(\left \lceil \alpha \right \rceil+2)} +\varepsilon }, \quad t\ge T_{\varepsilon}(\omega ),\quad \mathbb{P} -a.s.,
\end{split}
	\end{eqnarray*}
where  the random time $T_{\varepsilon}(\omega )$ is $\mathbb{P}-$a.s. finite;
\end{enumerate}
\item For any  $\varphi\in \mathcal{H}$  and $F\in C_{\mathcal{K} }(\mathcal{H} ) $ with $\int _{\mathcal{H}}F(\phi )\mu ^{*}(\text{d}\phi )=0$, let $\zeta =(\int _{\mathcal{H}}\mathcal{S}  [F(\phi )]\mu ^{*}(\text{d}\phi ))^{\frac{1}{2} }\in [0,\infty )$, then
\begin{enumerate}[(\textbf{a})]
		\item When $\zeta>0$, for $\epsilon  \in [0,\frac{1}{5} )$, there exists an increasing function $\mathcal{G} _{\epsilon  }:\mathbb{R}_{+}\times \mathbb{R} _{+}\to \mathbb{R} _{+}$ such that
\begin{eqnarray*}
\begin{aligned}
 \sup_{y\in \mathbb{R} } \left | \mathbb{P} (\frac{1}{\sqrt{t} }\int_{0}^{t}F(u_{r}(\varphi ))\text{d}r\le y  )-\Pi _{\zeta } (y)\right | \le \mathcal{G}_{\epsilon  }(\left \| F \right \| _{\mathcal{K} },\left \| \varphi  \right \|_{\mathcal{H} } )t^{-\frac{1}{5}+\epsilon   };
\end{aligned}
\end{eqnarray*}
\end{enumerate}
\begin{enumerate}[(\textbf{b})]
        \item When $\zeta=0$, there exists an increasing function $\mathcal{G} :\mathbb{R} _{+}\times \mathbb{R} _{+}\to \mathbb{R} _{+}$ such that
\begin{eqnarray*}
\begin{aligned}
 \sup_{y\in \mathbb{R} }[ (\left | y \right |\wedge 1 ) \left | \mathbb{P} (\frac{1}{\sqrt{t} }\int_{0}^{t}F(u_{r}(\varphi ))\text{d}r\le y  )-\Pi _{0 } (y)\right | ]\le \mathcal{G}(\left \| F \right \| _{\mathcal{K} },\left \| \varphi  \right \|_{\mathcal{H} } )t^{-\frac{1}{4} }.
\end{aligned}
\end{eqnarray*}
	\end{enumerate}
	\end{enumerate}
	\end{enumerate}
~\\
\\\textbf{proof}} The statements follow from Theorems 3.2, 4.5, 5.2, 6.1 and 6.6. Next, we need to show that the conditions \textbf{(A1)}$-$\textbf{(A6)} hold  for the above system.

Let $V=L^{q}(\Lambda  )$, $U=W^{1,2}_{0}(\Lambda )$ and
\begin{eqnarray*}
\begin{split}
		_{V ^{\ast}} \langle A(u), v\rangle _{V}:=-\int _{\Lambda }u(x)\left | u(x) \right |^{q-2}v(x)\text{d}x-a\int _{\Lambda }u(x)v(x)\text{d}x,
\end{split}
	\end{eqnarray*}
for $u,v\in V$, which implies that $V\subset  U= U^{*}\subset V^{*}$. Then we obtain:
\\(\textbf{A1}) First, by H\"older's inequality and Young's inequality, we obtain
\begin{eqnarray*}
\begin{split}
		\left |_{V ^{\ast}} \langle A(u_{1}), u_{2}\rangle _{V}  \right |  &\le [\int _{\Lambda }\left | u_{1}(x) \right |^{q}\text{d}x]^{1-\frac{1}{q} } \cdot[\int _{\Lambda }\left | u_{2}(x) \right |^{q}\text{d}x]^{\frac{1}{q} } \\
&~~~+a_{1}[\int _{\Lambda }\left | u_{1}(x) \right |^{\frac{q}{q-1} }\text{d}x]^{1-\frac{1}{q} } \cdot[\int _{\Lambda }\left | u_{2}(x) \right |^{q}\text{d}x]^{\frac{1}{q} }\\
&\le[\left \| u_{1} \right \| _{V}^{q-1}+a_{1}M_{q}(\left \| u_{1} \right \| _{V}^{q-1}+\left | \Lambda  \right |^{\frac{q-1}{q} })]\left \| u_{2} \right \| _{V},
\end{split}
	\end{eqnarray*}
which implies
\begin{eqnarray*}
\begin{split}
		\left \| A(u) \right \| _{V^{*}}\le(1+a_{1}M_{q})\left \| u \right \| _{V}^{q-1}+a_{1}M_{q}\left | \Lambda  \right |^{\frac{q-1}{q} },
\end{split}
	\end{eqnarray*}
i.e., $\gamma _{1}=(1+a_{1}M_{q})$ and $M=a_{1}M_{q}\left | \Lambda  \right |^{\frac{q-1}{q} }$.
\\(\textbf{A2})For all $u \in V $, we have by $a_{1}>0$
\begin{eqnarray*}
\begin{split}
		\left |_{V ^{\ast}} \langle A(u), u\rangle _{V}  \right |  \le-\left \| u \right \|^{q} _{V},
\end{split}
	\end{eqnarray*}
which implies $\gamma _{2}=0$, $p=q$ and $\gamma _{3}=1$.
\\(\textbf{A3}) For all $u_{1},u_{2} \in V $, we obtain
\begin{eqnarray*}
\begin{split}
		_{V ^{\ast}} \langle A(u_{1})-A(u_{2}), u_{1}-u_{2}\rangle _{V}    \le-a_{1}\left \| u_{1}-u_{2} \right \|^{2} _{U},
\end{split}
	\end{eqnarray*}
which implies $\lambda =-a_{1}$.
\\(\textbf{A4}) By the definition of $_{V ^{\ast}} \langle \cdot , \cdot \rangle _{V} $,  $\theta  \in \mathbb{R} \to_{V ^{\ast}}\langle A(u_{1}+\theta u_{2}),u_{3}\rangle _{V}$ is continuous  for all $u_{1}, u_{2}, u_{3} \in V $.
\\(\textbf{A6}) Let $\widetilde{U} =L^{2}(\Lambda ) $ and  $\Delta$ be the Laplace operator on $L^{2}(\Lambda ) $ with the Dirichlet boundary condition. Define $A_{n}=-\Delta(I-\frac{\Delta }{n} )^{-1}$, then $\textbf{A6}$ holds(see \cite{R3} for details).
~\\

For example, let
\begin{eqnarray*}
\begin{split}
		f(t, u_{t})= a_{2}(\sin t+\cos 3t )[\cos u(t,x)+\int_{-\tau   }^{0}u(t+\theta,x )\text{d}\theta ]  ,
\end{split}
	\end{eqnarray*}
and
\begin{eqnarray*}
\begin{split}
		g(t, u_{t})= a_{3}\sin (\frac{\sqrt{2}}{2+\cos t+\cos \sqrt{2}t) } )\int_{-\tau  }^{0}u(t+\theta,x )\text{d}\theta,
\end{split}
	\end{eqnarray*}
where $a_{2},a_{3}\in \mathbb{R} $. For any $u^{1}_{t} , u^{2}_{t} \in \mathcal{H} $, we have
\begin{eqnarray*}
		\left \langle f(t,u^{1}_{t} )-f(t,u^{2}_{t}),u^{1}(t)-u^{2}(t)\right \rangle_{U}  \le 2\left | a_{2} \right | \cdot  [\left \| u^{1}(t)-u^{2}(t) \right \|_{U} ^{2} +\int_{-\tau  }^{0}\left \| u^{1}(t+\theta )-u^{2}(t+\theta ) \right \|^{2}_{U}  \text{d}\theta ],
	\end{eqnarray*}
\begin{eqnarray*}
		\left \| f(t,u^{1}_{t} )-f(t,u^{2}_{t} ) \right \|_{U} \le 4\left | a_{2} \right | \cdot \left \| u^{1}_{t} -u^{2}_{t} \right \| _{\mathcal{H} },
	\end{eqnarray*}
\begin{eqnarray*}
		\left \| g(t,u^{1}_{t} )-g(t,u^{2}_{t} ) \right \| _{\mathscr{L}(K,U)}^{2}  \le \left | a_{3} \right |\cdot \int_{-\tau  }^{0}\left \| u^{1}(t+\theta )-u^{2}(t+\theta ) \right \|^{2}_{U}  \text{d}\theta ,
	\end{eqnarray*}
which implies $\eta _{1}=-2\left | a_{2} \right | $, $\eta _{2}=2\left | a_{2} \right |$, $L_{0}=4\left | a_{2} \right |$ and $\eta _{3}=\left | a_{3} \right |$. Thus, when $a_{1}-4\left | a_{2} \right |-\frac{145\left | a_{3} \right | }{2}  >0$, there exists a unique Levitan almost periodic measure $\mu _{t}\in \mathcal{R}_{2}$ in the sense of Wasserstein metric for  system \eqref{g31}. In addition,  let
		$$f(t, u_{t})= a_{2}[\cos u(t,x)+\int_{-\tau   }^{0}u(t+\theta,x )\text{d}\theta ],\quad
		g(t, u_{t})= a_{3}\int_{-\tau  }^{0}u(t+\theta,x )\text{d}\theta,$$
then $\eta _{1}=-\left | a_{2} \right | $, $\eta _{2}=\left | a_{2} \right |$, $L_{0}=2\left | a_{2} \right |$ and $\eta _{3}=\left | a_{3} \right |$. Thus, when $a_{1}-2\left | a_{2} \right |-\frac{145\left | a_{3} \right | }{2}  >0$, the conclusions (i)-(iii) of Theorem 7.1 are also hold.

\section*{Appendix \uppercase\expandafter{\romannumeral1}: The specific proof of Lemma 4.1:}
For any $t\ge s$, let $\varsigma_{0} (t)$ be a stopping time as defined by
\begin{eqnarray*}
\begin{split}
		\varsigma_{0} (t)=\inf\{r\ge t:\left \| u(r)\right \|>R_{0}\} ,
\end{split}
	\end{eqnarray*}
and we set $\varsigma_{0} (t)=+\infty $ if $\inf\{r\ge t:\left \| u(t) \right \|>R_{0}\}=\emptyset $. For any $T \ge 0$, by $\textbf{(A2)}$, $\textbf{(A5)}$ and \eqref{51},\eqref{26}, applying the It$\hat{\text{o}} $ formula to $ \left \| u((t+T)\wedge\varsigma_{0} (t) ) \right \| ^{2} $   yields
\begin{align}\label{r30}
\begin{split}
		&\mathbb{E} \left \| u((t+T)\wedge\varsigma_{0} (t) ) \right \| ^{2}\\&=\mathbb{E} \left \| u(t) \right \| ^{2}+\mathbb{E} \int_{t}^{(t+T)\wedge\varsigma_{0} (t)}[2_{V ^{\ast}} \langle A(r,u(r)),u(r))\rangle _{V}+2\left \langle f(r,u_{r}),u(r) \right \rangle \\&~~~+\left \| g(r,u_{r}) \right \| ^{2}_{\mathscr{L}(K,U)}    ]\text{d}r\\&\le L_{1}+L_{2}\left \| \varphi  \right \| _{\mathcal{H}}^{2 }+[-(2\eta _{1}-\epsilon_{1})+2\lambda ]\int_{t}^{(t+T)\wedge\varsigma_{0} (t)}\mathbb{E} \left \| u(r) \right \| ^{2}\text{d}r\\&~~~+(2\eta _{2}+\frac{\eta_{3}}{1-\epsilon_{2}} )\int_{t}^{(t+T)\wedge\varsigma_{0} (t)}\int_{-\tau }^{0} \mathbb{E} \left \| u(r+\theta ) \right \| ^{2}\pi  (\text{d}\theta )\text{d}r+(\frac{1}{\epsilon_{1}}+\frac{1}{\epsilon_{2}} )M^{2}T\\&\le L_{1}+L_{2}\left \| \varphi  \right \| _{\mathcal{H}}^{2 }+[-(2\eta _{1}-\epsilon_{1})+2\lambda +2\eta _{2}+\frac{\eta_{3}}{1-\epsilon_{2}} ]\int_{t}^{(t+T)\wedge\varsigma_{0} (t)}\mathbb{E} \left \| u(r) \right \| ^{2}\text{d}r\\&~~~+(2\eta _{2}+\frac{\eta_{3}}{1-\epsilon_{2}} )\int_{t-\tau }^{t} \mathbb{E} \left \| u(r+\theta ) \right \| ^{2}\pi  (\text{d}\theta )\text{d}r+(\frac{1}{\epsilon _{1}}+\frac{1}{\epsilon _{2}} )M^{2}T\\&\le L_{1}+L_{2}\left \| \varphi  \right \| _{\mathcal{H}}^{2 }-\tau(2\eta _{2}+\frac{\eta_{3}}{1-\epsilon_{2}} )\frac{G_{1}+G_{4}}{G_{2}+G_{5}}+(\frac{1}{\epsilon _{1}}+\frac{1}{\epsilon_{2}} )M^{2}T\\&:=\overline{L} (\varphi ,T).
\end{split}\tag{\uppercase\expandafter{\romannumeral1}.1}
	\end{align}
By \eqref{r30}, we have
\begin{align}\label{r31}
\begin{split}
		\mathbb{E} [\left \| u(\varsigma_{0} (t)) \right \|^{2}\textbf{\emph{1}} _{\varsigma_{0} (t)<t+T} ]\le \overline{L}  ,
\end{split}\tag{\uppercase\expandafter{\romannumeral1}.2}
	\end{align}
which implies
\begin{align}\label{r32}
\begin{split}
		\mathbb{P} (\omega :\varsigma_{0} (t)<t+T)\le \frac{\overline{L}  }{R_{0}^{2}} .
\end{split}\tag{\uppercase\expandafter{\romannumeral1}.3}
	\end{align}
Hence for any $t\ge s$, there exists a positive constant $R' _{0}$ independent of $\epsilon\in[0, 1]$ such that the solution $u(t;s,\varphi )$ of \eqref{rrr1} satisfies
\begin{eqnarray*}
		\mathbb{P} (\omega :\{\sup_{r\in [t,t+T]}\left \| u(r;s,\varphi ) \right \|>R  \})<\epsilon, \quad t\ge s,\quad  R\ge R' _{0}.
	\end{eqnarray*}
This completes the proof. \quad $\Box $

\section*{Appendix \uppercase\expandafter{\romannumeral2}: The specific proof of Lemma 4.2:}
Consider the  system \eqref{rrr1},
\begin{eqnarray*}
		\text{d}u(t)=[A(t,u(t))+f(t,u_{t})dt+g(t,u_{t})]\text{d}W(t), \quad t\ge s
	\end{eqnarray*}
with the initial data $u_{s} =\varphi \in \mathcal{H}$. We know that under $\textbf{(A1)}$-$\textbf{(A5)}$, \eqref{rrr1} admits a unique solution $u(t;s,\varphi):=u(t)$ and note that
\begin{eqnarray*}
\begin{split}
		u(t;s,\varphi)=\varphi(0)+\int_{s}^{t}[A(r,u(r;s,\varphi ))+ f(r,u_{r}(s,\varphi ))]\text{d}r+\int_{s}^{t} g(r,u_{r}(s,\varphi ))\text{d}W(r).
\end{split}
	\end{eqnarray*}
For any $\epsilon > 0$, letting $T=2\tau $ of \textbf{Lemma 4.1},  it follows from \textbf{Proposition 1} and the definition of the norm of
$\widetilde{U} $ that there exists $R_{2} > 0$ such that for any $t\ge s$
\begin{align}\label{r21}
\begin{split}
		\mathbb{P} (\omega :\sup_{r\in [t,t+\tau ]}\left \| u_{t} \right \|_{\tilde{\mathcal{H} }}>R_{2}  )<\epsilon.
\end{split}\tag{\uppercase\expandafter{\romannumeral2}.1}
	\end{align}

For any $t\ge s$, let $\varsigma_{1} (t)$ be a stopping time as defined by
\begin{align}\label{r22}
\begin{split}
		\varsigma_{1} (t)=\inf\{r\ge t:\left \| u_{t} \right \|_{\tilde{\mathcal{H} }}>R_{2}\} ,
\end{split}\tag{\uppercase\expandafter{\romannumeral2}.2}
	\end{align}
and we set $\varsigma_{1} (t)=+\infty $ if $\inf\{r\ge t:\left \| u_{t} \right \|_{\tilde{\mathcal{H} }}>R_{2}\}=\emptyset $.
For any $t\ge s$, $\delta \in(0,\tau  )$ and $l>1$, by $\textbf{(A2)}$, $\textbf{(A5)}$ and Burkholder–Davis–Gundy inequality, we have
\begin{eqnarray*}
\begin{split}
		&\mathbb{E} (\sup_{\theta \in [t,t+\delta ]}\left \| u(\theta\wedge \varsigma_{1} (t) )-u(t) \right \|^{2l}  )\\&\le 3^{2l-1}\mathbb{E} (\sup_{\theta \in [t,t+\delta ]}\left \| \int_{t}^{\theta\wedge \varsigma_{1} (t) }A(r,u(r;s,\varphi ))\text{d}r  \right \|^{2l}  )+3^{2l-1}\mathbb{E} (\sup_{\theta \in [t,t+\delta ]}\left \| \int_{t}^{\theta\wedge \varsigma_{1} (t) }f(r,u_{r}(s,\varphi ))\text{d}r  \right \|^{2l}  )\\&~~~+3^{2l-1}\mathbb{E} (\sup_{\theta \in [t,t+\delta ]}\left \| \int_{t}^{\theta\wedge \varsigma_{1} (t) }g(r,u_{r}(s,\varphi ))\text{d}W(r)  \right \|^{2l}  )\\&\le 3^{2l-1}[(\gamma _{1}R_{1}^{p-1}+ M)^{2l}+(M+L_{0}R_{1})^{2l}]\delta ^{2l}+3C_{l}\mathbb{E} ( \int_{t}^{(t+\delta)\wedge \varsigma_{1} (t) }\left \|g(r,u_{r}(s,\varphi ))\right \|^{2}\text{d}r )^{l}\\&\le G_{7}\delta^{l}.
\end{split}
	\end{eqnarray*}
Let $s\to-\infty $ such that for any $t\in \mathbb{R} $
\begin{align}\label{r19}
\begin{split}
	\mathbb{E} (\sup_{\theta \in [t,t+\delta ]}\left \| \mathcal{U}  (\theta\wedge \varsigma_{1} (t) )-\mathcal{U}  (t) \right \|^{2}  )\le	G_{7}\delta^{l}.
\end{split}\tag{\uppercase\expandafter{\romannumeral2}.3}
	\end{align}
 By \eqref{r21} and \eqref{r22}, we obtain $\mathbb{P} (\omega :\varsigma_{1} (t-\tau)<t)<\epsilon$. Hence let $\delta =1\wedge \tau \wedge (\frac{\epsilon\kappa ^{2l}}{\tau 9^{l}} )^{\frac{1}{l-1} } $, then we have
\begin{eqnarray*}
\begin{split}
&\mathbb{P} (\omega :\{\sup_{\theta _{1},\theta _{2}\in [-\tau ,0],\left | \theta _{1}-\theta _{2} \right |<\delta  }\left \| \mathcal{U}  _{t}(\theta _{1})-\mathcal{U}  _{t}(\theta _{2}) \right \|\ge \kappa \}  )\\&\le \mathbb{P}(\omega :\varsigma_{1} (t-\tau)<t)+\mathbb{P}(\omega :\{\varsigma_{1} (t-\tau)>t,\sup_{\theta _{1}\in [-\tau ,0],\theta _{2} \in[\theta _{1},(\theta _{1}+\delta )\wedge 0]  }\left \| \mathcal{U}  (t+\theta _{1})-\mathcal{U} (t+\theta _{2}) \right \|\ge \kappa \}  )\\&\le\varepsilon+\mathbb{P}(\omega :\{\varsigma_{1} (t-\tau)>t,\max_{k\in [0,\left [ \frac{\tau }{\delta }  \right ] ],k\in \mathbb{N}^{+} } \sup_{r\in [t-(k+1)\delta\wedge \tau  ,t-k\delta] }\left \| \mathcal{U}  (r)-\mathcal{U} (t-(k+1)\delta\wedge \tau) \right \|\ge \frac{\kappa}{3 } \}  )\\&\le\varepsilon+\sum_{k=0}^{\left [ \frac{\tau }{\delta }  \right ] } \mathbb{P}(\omega :\{\varsigma _{1} (t-(k+1)\delta\wedge \tau)>t,\sup_{r\in [t-(k+1)\delta\wedge \tau  ,t-k\delta] }\left \| \mathcal{U}  (r)-\mathcal{U} (t-(k+1)\delta\wedge \tau) \right \|\ge \frac{\kappa}{3 } \}  )\\&\le\varepsilon+\sum_{k=0}^{\left [ \frac{\tau }{\delta }  \right ] } \mathbb{P}(\omega :\{\sup_{r\in [t-(k+1)\delta\wedge \tau  ,t-k\delta] }\left \|\mathcal{U} (r\wedge \varsigma _{1} (t-(k+1)\delta\wedge \tau))-\mathcal{U} (t-(k+1)\delta\wedge \tau) \right \|\ge \frac{\kappa}{3 } \}  )\\&\le \varepsilon+(1+\frac{\tau }{\delta } )\frac{3^{2l}}{\kappa ^{2l}} G_{7}\delta ^{l}\\&\le G_{8}\varepsilon.
\end{split}
	\end{eqnarray*}
This completes the proof. \quad $\Box $

\section*{Acknowledgments}
 The first author (S. Lu)  supported by Graduate Innovation Fund of Jilin University. The second author (X. Yang) was supported by  National Natural Science Foundation of China (12071175, 12371191). The third author (Y. Li) was supported by  National Basic Research Program of China (2013CB834100), National Natural Science Foundation of China (12071175, 11171132 and 11571065), Project of Science and Technology Development of Jilin Province (2017C028-1 and 20190201302JC) and Natural Science Foundation of Jilin Province (20200201253JC).
\section*{Data availability}
No data was used for the research described in the article

\end{document}